\documentclass[11pt, reqno]{amsart}

\usepackage{amssymb}
\usepackage{amsmath}
\usepackage{amsxtra}
\usepackage{amscd}
\usepackage{color}
\usepackage{cases}

\newtheorem{theorem}{Theorem}
\newtheorem{corollary}{Corollary}
\newtheorem{lemma}{Lemma}
\newtheorem{proposition}{Proposition}
\newtheorem{definition}{Definition}
\newtheorem{remark}{Remark}

\begin{document}

\title{The Yokonuma--Temperley--Lieb Algebra}

\author{D. Goundaroulis}
\address{Department of Mathematics,
National Technical University of Athens,
Zografou campus, GR-157 80 Athens, Greece.}
\email{dgound@mail.ntua.gr}
\urladdr{users.ntua.gr/dgound}

\author{J. Juyumaya}
\address{Instituto de Matem\'{a}ticas, Universidad de Valpara\'{i}so,
Gran Breta\~{n}a 1091, Valpara\'{i}so, Chile}
\email{juyumaya@uvach.cl}

\author{A. Kontogeorgis}
\address{Department of Mathematics, University of Athens\\
Panepistimioupolis, 15784 Athens, Greece
}
\email{kontogar@math.uoa.gr}

\author{S. Lambropoulou}
\address{ Department of Mathematics,
National Technical University of Athens,
Zografou campus, GR-157 80 Athens, Greece.}
\email{sofia@math.ntua.gr}
\urladdr{http://www.math.ntua.gr/~sofia}


\thanks{This research has been co-financed by the European Union (European
Social Fund - ESF) and Greek national funds through the Operational Program
``Education and Lifelong Learning" of the National Strategic Reference
Framework (NSRF) - Research Funding Program: THALIS. Moreover, the second author was partially supported by the National Technical U. Athens and DIUV
 No.01-2011}


\begin{abstract}
In this paper we introduce the Yokonuma--Temperley--Lieb algebra as a quotient of the Yokonuma--Hecke algebra over a two-sided ideal generated by an expression analogous to the one of the classical Temperley-Lieb algebra. The main theorem provides necessary and sufficient conditions for the Markov trace defined on the Yokonuma--Hecke algebra to pass through to the quotient algebra, leading to a sequence of knot invariants which coincide with the Jones polynomial.
\end {abstract}

\date{}
\maketitle

\section*{Introduction} \label{s:intro}
The Temperley--Lieb algebra appeared originally  in Statistical Mechanics and is 
important in several areas of Mathematics. In his seminal work V.F.R. Jones \cite{JonesIndex} constructed a Markov trace  on the Temperley--Lieb algebra, leading  to unexpected applications in knot theory as well as to a  fertile  interaction between Knot theory and Representation theory. In algebraic terms, the Temperley--Lieb algebra can be defined as a quotient of the Iwahori--Hecke algebra.

\smallbreak
In \cite{Juyu} the Yokonuma--Hecke algebra ${\rm Y}_{d,n}(u)$ (defined originally in \cite{Yo}) has been defined as a quotient of the modular framed braid group $\mathcal{F}_{d,n}$, which comprises framed braids with framings modulo $d$, over a quadratic relation (Eq.~\ref{sigmaquadr}) involving the framing generators $t_i$ by means of certain weighted idempotents $e_{i}$ (Eq.~\ref{edi}). Setting $d = 1$, the algebra ${\rm Y}_{1,n}(u$) coincides with the Iwahori--Hecke algebra ${\rm H}_n(u)$. The Yokonuma--Hecke  algebras have been studied in \cite{Yo,Juyu,JL,th,ChLo}. Further, in \cite{Juyu} the second author found an inductive linear basis for the algebras ${\rm Y}_{d,n}(u)$ and  constructed a unique Markov trace  {\rm tr} on these algebras depending on parameters $z, x_1, \ldots, x_{d-1}$. Aiming to extracting framed link invariants from {\rm tr}, as it turned out in \cite{JL2}, ${\rm tr}$ does not re-scale directly according to the framed braid equivalence, leading to conditions that have to be imposed on the trace parameters $x_1, \ldots, x_{d-1}$; namely, they had to satisfy a non-linear system of equations, the {\it ${\rm E}$--system} (Eq.~\ref{Esystem}). The $x_i$'s being $d^{th}$ roots of unity is one obvious solution. G\'{e}rardin found in \cite[Appendix]{JL2} the full set of solutions of the ${\rm E}$--system. Given now any solution of the ${\rm E}$--system, 2--variable isotopy invariants for framed, classical and singular links  were constructed  in \cite{JL2,JL4,JL3} respectively, which are studied further in \cite{ChLa,CJJKL}.

\smallbreak
In this paper we define a Temperley--Lieb analogue of the Yokonuma--Hecke algebra, {\it the Yokonuma--Temperley--Lieb algebra} ${\rm YTL}_{d,n}(u)$, as a quotient of the Yokonuma--Hecke algebra over a two-sided ideal  $I$ (Eq.~\ref{ideal} and Definition~\ref{ytldef}), analogous to the classical case. For $d=1$ the algebra ${\rm YTL}_{1,n}(u)$ coincides with the Temperley--Lieb algebra. We first show that $I$ is a principal ideal (Corollary~\ref{coroldef}) and we give a presentation for ${\rm YTL}_{d,n}(u)$ with non-invertible generators, analogous to the classical case (Proposition \ref{noninvert}). We then give a spanning set $\Sigma_{d,n}$ for ${\rm YTL}_{d,n}(u)$, where each word in $\Sigma_{d,n}$ contains the highest and lowest index braiding generator exactly once (Proposition~\ref{ytlspan}). Moreover, any word in $\Sigma_{d,n}$ inherits the splitting property from ${\rm Y}_{d,n}(u)$, that is, it splits into the framing part  and the braiding part. We also present the results of Chlouveraki and Pouchin \cite{ChPou} on the dimension (Proposition~\ref{ytldim}) and a linear basis for ${\rm YTL}_{d,n}(u)$ (Theorem~\ref{chpouthm}). From the spanning set $\Sigma_{d,n}$, they extracted an explicit basis for ${\rm YTL}_{d,n}(u)$ by describing a set of linear dependence relations among the framing parts for each fixed element in the braiding part. Finally, using the dimension results of \cite{ChPou} we find a basis for ${\rm YTL}_{2,3}(u)$ different than the basis in \cite{ChPou}\\

Next, we seek conditions such that the trace ${\rm tr}$, defined on the algebras ${\rm Y}_{d,n}(u)$, passes to the quotient algebras ${\rm YTL}_{d,n}(u)$. More precisely, we compute first the values of the trace parameter $z$ that annihilate the generator of the defining ideal $I$, which are the roots of a quadratic equation (Eq.~\ref{trg120}). Then we annihilate the traces of all elements of ${\rm Y}_{d,n}(u)$ that lie in $I$ and so we end up with a system $(\Sigma)$ of quadratic equations in $z$ (Eqs.~\ref{redeq1}--\ref{redeq4}).  If we demand that $(\Sigma)$ has both roots of Eq.~\ref{trg120} as common solutions, which is essential for discussing link invariants, we end up with necessary conditions for the trace ${\rm tr}$ to pass to the quotient algebras ${\rm YTL}_{d,n}(u)$ (Theorem~\ref{generalcase}). More precisely, Theorem~\ref{generalcase} states that the trace ${\rm tr}$ passes to the quotient algebra ${\rm YTL}_{d,n}(u)$ if the trace parameters are $d^{th}$ roots of unity  $x_1 ,\ldots , x_{d-1}$ and $z=-\frac{1}{u+1}$ and $z=-1$. Note that these two values for $z$ are precisely the ones that Jones computed such that the Ocneanu trace on ${\rm H}_n(u)$ passes to the quotient, the Temperley--Lieb algebra ${\rm TL}_n(u)$. If we also let $(\Sigma)$  to have one common solution for $z$ we obtain the necessary and sufficient conditions for the trace ${\rm tr}$ to pass through to the quotient algebras ${\rm YTL}_{d,n}(u)$ (Theorem~\ref{ytlthmgen}). More precisely, Theorem~\ref{ytlthmgen} states that the trace {\rm tr} passes to the quotient algebras ${\rm YTL}_{d,n}(u)$ if and only if either the conditions of Theorem~\ref{generalcase} are satisfied or the trace parameters $x_1 , \ldots , x_{d-1}$ comprise a solution of the ${\rm E}$--system (other than $d^{th}$ roots of unity) and $z= - \frac{1}{2}$. This is our main result.\\

In \cite{ChLa} it is shown that if the trace parameters $x_1 , \ldots , x_{d-1}$ are $d^{th}$ roots of unity, then the classical link invariants derived from the algebra ${\rm Y}_{d,n}(u)$ coincide with the 2--variable Jones or HOMFLYPT polynomial. Using Theorem~\ref{ytlthmgen} and the results in \cite{ChLa}, we obtain from the invariants for framed and classical links in \cite{JL2,JL4} related to ${\rm Y}_{d,n}(u)$  1--variable framed and classical link invariants through the algebras ${\rm YTL}_{d,n}(u)$ (Definition~\ref{invsdef}). As we show, these invariants coincide with the Jones polynomial for the case of classical links and they are framed analogues of the Jones polynomial for the case of framed links (Corollary~\ref{invcor}).\\

The paper is organized as follows: In Section~1 we recall the definition and basic properties of the classical Temperley--Lieb algebra and the Yokonuma--Hecke algebra. In Section~2 we define the Yokonuma--Temperley--Lieb algebra as a quotient of the Yokonuma--Hecke algebra over a two-sided ideal (Eq.~\ref{ideal} and Definition~\ref{ytldef}), which we show that is a principal ideal (Corollary~\ref{coroldef}).  Finally, we give a presentation for ${\rm YTL}_{d,n}(u)$ with non-invertible generators (Proposition \ref{noninvert}).  In Section~3 we present a spanning set for ${\rm YTL}_{d,n}(u)$ and the results of Chlouveraki and Pouchin \cite{ChPou} on the dimension and a linear basis for ${\rm YTL}_{d,n}(u)$. Then we give a basis for ${\rm YTL}_{2,3}(u)$. Section~4 focuses on the necessary and sufficient conditions under which the trace ${\rm tr}$ on ${\rm Y}_{d,n}(u)$ passes to the quotient algebra ${\rm YTL}_{d,n}(u)$ (Theorems~\ref{generalcase} and \ref{ytlthmgen}). Finally, in Section~5 we discuss the invariants for classical and framed links that can be constructed through the trace ${\rm tr}$ and we recover the Jones polynomial (Corollary~\ref{invcor}).

\section{Preliminaries}
\subsection{\it Notations} Throughout the paper we shall fix the following notation. By the term algebra we mean an associative unital (with unity 1) algebra over the field $K:= \mathbb{C}(u)$, where $u$ is an indeterminate. The following two possitive integers are also fixed: $d$ and $n$. 

As usual we denote by $B_n$ the braid group  on $n$ strands, that is the group generated 
by the elementary braids $\sigma_1,\ldots ,\sigma_{n-1}$, where
$\sigma_i$ is the positive crossing between the $i^{th}$ and the
$(i+1)^{st}$ strand, satisfying the well-known braid relations: $\sigma_i \sigma_{i+1} \sigma_i = \sigma_{i+1} \sigma_i \sigma_{i+1}$ and $\sigma_i \sigma_j = \sigma_i \sigma_j$ for $|i-j|>1$. 

We denote $S_n$ the symmetric group on $n$ symbols. Let $s_i$ be the elementary transposition $(i,i+1)$. We denote by $l$ the length function on $S_n$ with respect to the $s_i$'s.

Denote by $C_d = \langle t \,| \,t^d=1 \rangle$ the cyclic group of order $d$. Let $t_i = (1 ,\ldots ,t ,1, \ldots , 1)\in C_d^n$, where $t$ is in the $i^{th}$ position. 

Finally, we denote $C_{d,n }:= C_d^n \rtimes S_n$, where the action is defined by permutation on the indices of the $t_i$'s, namely: $s_i t_j = t_{s_i(j)}s_i$.

\subsection{\it{The  Temperley--Lieb algebra}}\label{clasTL}

Originally, the Temperley--Lieb algebra, over $\mathbb{C}$, was defined by generators $ 1, f_1, \ldots , f_{n-1} $ subject to the following relations:
$$
\begin{array}{ccl}
f_{i}^{2} & = &  f_i \\
f_if_{j}f_i &  = &  \tau f_i, \quad |i-j| = 1\\
f_if_j &  = & f_j f_i,  \quad |i-j| > 1
\end{array}
$$
where $\tau$ is a non--zero complex number (see  \cite{Goodman},\cite{Jones},\cite{JonesIndex}). 
The generators $f_i$ are non-invertible; 
 one can define the Temperley--Lieb algebra  with the following invertible generators (see \cite{Jones}):
 \begin{equation} \label{transf} 
h_i := (q+1)f_i-1
\end{equation}
where 
 $q$ is defined via the relation $\tau^{-1} = 2+ q +q^{-1}$. The Temperley algebra ${\rm TL}_n(u)$, over $K$, is defined by generators $h_1, \ldots ,h_{n-1}$ under the relations:

\begin{align}
h_ih_jh_i &= h_jh_ih_j, \quad |i-j| =1 \label{TL1}\\
h_ih_j &= h_jh_i, \quad |i-j|> 1\label{TL2} \\
h_i^2 &= (u-1)h_i + u  \label{TL3}\\
h_ih_{j}h_i + h_{j}h_i +& h_i h_{j} + h_i  + h_{j} + 1= 0, \quad |i-j|=1. \label{TL4}
\end{align}

\smallbreak

Note that relations (\ref{TL4}) are symmetric with respect to the indices $i,j$, so relations (\ref{TL1}) follow from relations (\ref{TL4}).  Relations (\ref{TL1})--(\ref{TL3}) are the well--known defining relations of the Iwahori--Hecke algebra ${\rm H}_n(u)$.  Therefore, ${\rm TL}_n(u)$  can be considered as a quotient of ${\rm H}_{n}(u)$ over the two--sided ideal generated by relations (\ref{TL4}). It turns out that the set:
$$
\left \{ \left ( h_{j_1}h_{j_1 -1} \ldots h_{j_1 - k_1} \right ) \left( h_{j_2}h_{j_2 -1} \ldots h_{j_2 - k_2} \right )\ldots \left ( h_{j_p}h_{j_p -1} \ldots h_{j_p - k_p} \right ) \right \} $$
where $1 \leq j_1 < j_2 < \ldots <j_p \leq n-1 $ and $ 1 \leq j_1 - k_1 < j_2 - k_2 < \ldots < j_p - k_p $,
furnishes a linear basis for ${\rm TL}_n(u)$ and the dimension of ${\rm TL}_n(u)$ is equal to the $n^{th}$ Catalan number $c_n= \frac{1}{n+1}$~$2n\choose{n}$ \cite{Jones, JonesIndex}. Recall finally, that in \cite{Homfly}, Ocneanu constructed a unique Markov trace on the algebras ${\rm H}_n (u)$:
\begin{theorem}[Ocneanu] For any $\zeta \in K^{{\times}}$ there exists a linear trace ${\rm \tau}$ on $\cup_{n=1}^{\infty} {\rm H}_n(u)$ uniquely defined by the inductive rules:
\begin{enumerate}
\item ${\rm \tau} (a b) = {\rm \tau} (ba), \quad a,b \in {\rm H}_n(u)$
\item ${\rm \tau} (1) = 1$
\item ${\rm \tau} (a g_n ) = \zeta \, {\rm \tau}(a), \quad a\in {\rm H}_n(u)$.
\end{enumerate}
\end{theorem}
Jones' methods for redefining his Markov trace on the Temperley--Lieb algebra as factoring of the Ocneanu trace on the Iwahori--Hecke algebra \cite{Jones} tells us that the least requirement is that the Ocneanu trace respects the defining relations (\ref{TL4}). This requirement implies:
\begin{equation}\label{jonval}
\zeta = - \frac{1}{u+1} \quad \mbox{and} \quad \zeta=-1.
\end{equation}
The Ocneanu trace is used in \cite{Jones} for constructing the HOMFLYPT polynomial invariant for classical knots and links. Then, by specializing $\zeta$ to $- \frac{1}{u+1}$ the Jones polynomial was recovered.
\subsection{\it{The Yokonuma--Hecke algebra}}
The group
${\mathbb Z}^n$ is generated by the  ``framing generators'' $t_1, \ldots , t_n$, the standard multiplicative generators of $\mathbb{Z}^n$. In this notation an element $a=(a_1, \ldots,a_n) \in {\mathbb Z}^n$ in the additive notation can be expressed as $ t_1^{a_1}\ldots t_{n}^{a_n}$. The {\it framed braid group} on $n$ strands is then defined as:
\begin{equation*}
{\mathcal F}_{n} = {\mathbb Z}^n \rtimes  B_n
\end{equation*}
where the action of $B_n$ on ${\mathbb Z}^n$ is given by the permutation induced by a braid on the indices:
\begin{equation}\label{action}
\sigma_it_j=t_{s_i(j)}\sigma_i.
\end{equation}
In particular, $\sigma_i t_i = t_{i+1}\sigma_i$ and $t_{i+1}\sigma_i = \sigma_i t_i$. A word $w$ in ${\mathcal F}_{n}$ has thus the ``splitting property'', that is, it splits into the  ``framing'' part and the ``braiding''  part:
$$w = t_1^{a_1}\ldots t_n^{a_n} \, \sigma$$  where $\sigma \in B_n$ and $a_i \in \mathbb{Z}$. So $w$ is a
classical braid with an integer attached to each
strand. Topologically, an element of ${\mathbb Z}^n$ is identified with
a framed identity braid on $n$ strands, while a classical braid in $B_n$ is
viewed as a framed braid with all framings 0. The multiplication in
${\mathcal F}_n$ is defined by placing one braid on top of the other and
collecting the total framing of each strand to the top.

For a fixed positive integer $d$, the {\it $d$-modular framed braid
group} on $n$ strands, ${\mathcal F}_{d,n}$, is defined as the
quotient of ${\mathcal F}_n$ over the {\it modular relations}:
\begin{equation}\label{modular}
t_i^d= 1 \quad (i =1, \ldots, n).
\end{equation}
Thus, ${\mathcal F}_{d,n}= C_d^n \rtimes  B_n$, where $C_d^n$ is isomorphic to $({\mathbb Z}/d{\mathbb Z})^n$ but with multiplicative notation.
Framed braids in ${\mathcal F}_{d,n}$  have framings modulo $d$. 

Passing now to the group algebra $\mathbb{C}{\mathcal F}_{d,n}$, we have the following elements  $e_{i} \in \mathbb{C}C_d^n$ (see \cite{JL} for diagrammatic interpretations), which are idempotents (cf. \cite[Lemma 4]{JL}):
\begin{equation}\label{edi}
e_{i} := \frac{1}{d} \sum_{s=0}^{d-1}t_i^s t_{i+1}^{-s},
\qquad i=1,\ldots , n-1.
\end{equation}
The definition of the idempotent $e_i$ can be generalized in the following way. For any indices $i,j$ and any $m \in \mathbb{Z}/ d \mathbb{Z}$, we define the following elements in ${\rm Y}_{d,n}(u)$:
\begin{equation}\label{eij}
e_{i,j} := \frac{1}{d} \sum_{s=0}^{d-1}t_i^s t_{j}^{-s},
\end{equation}
and:
\begin{equation}\label{eijm}
e_i^{(m)} := \frac{1}{d} \sum_{s=0}^{d-1} t_i^{m+s} t_{i+1}^{-s}.
\end{equation}
(notice that $e_i = e_{i, i+1} = e_i^{(0)}$).
The following lemma collects some of the relations among the $e_i$'s, the $t_i$'s and the $g_i$'s. These relations  will be used in the paper.
\smallbreak
\begin{lemma}\label{eipropop}
For the idempotents $e_i$  and for $1 \leq i , \,  j \leq n-1$ the following relations hold:
$$
\begin{array}{rcl}
t_j e_i &=& e_i t_j\\
e_{i+1}g_i &=& g_i e_{i,i+2}\\
e_ig_j &=& g_j e_i, \quad \mbox{for } j \neq i-1, i+1\\
e_jg_ig_j &=& g_ig_je_i \quad \mbox{for } |i-j|=1\\
e_ie_{i+1}&=& e_i e_{i, i+2}\\
e_{i}e_{i+1} &=& e_{i,i+2}e_{i+1}.
\end{array}
$$
\end{lemma}
\begin{proof}
All relations are immediate consequences of the definitions. The proofs for the first four relations can be found, for example, in \cite[Lemma 2.1]{JL3}. For the fifth relation we have:
\begin{eqnarray}
e_i e_{i+1}  &=& \frac{1}{d} \sum_{s=0}^{d-1} t_i^s t_{i+1}^{-s}\, \frac{1}{d} \sum_{m=0}^{d-1} t_{i+1}^m t_{i+2}^{-m}\nonumber\\
&=& \frac{1}{d^2} \sum_{s=0}^{d-1}\,\sum_{m=0}^{d-1} t_i^s t_{i+1}^{m-s} t_{i+2}^{-m}. \label{doublee}
\end{eqnarray}
Setting now $k=m-s$ we obtain:
\begin{eqnarray*}
(\ref{doublee}) &=& \frac{1}{d^2} \sum_{s=0}^{d-1}\,\sum_{k=0}^{d-1} t_i^s t_{i+1}^{k} t_{i+2}^{-k-s}\\
&=& \frac{1}{d} \sum_{s=0}^{d-1} t_i^s t_{i+2}^{-s} \,\frac{1}{d} \sum_{k=0}^{d-1} t_{i+1}^{k} t_{i+2}^{-k}\\
&=& e_{i, i+2} e_{i+1}.
\end{eqnarray*}
The sixth relation is proved in an analogous way.
\end{proof}
The {\it Yokonuma--Hecke algebra} ${\rm Y}_{d,n}(u)$ is defined \cite{Juyu, JL} as the quotient of the group algebra
$\mathbb{C} {\mathcal F}_{d,n}$ over the two-sided ideal  generated by the elements:
$$
\sigma_i^2 - 1 - (u-1) \, e_i - (u-1) \, e_i \, \sigma_i, \quad \mbox{for all } i,
$$
which give rise to the following quadratic relations in ${\rm Y}_{d,n}(u)$: 
\begin{equation}\label{sigmaquadr}
g_i^2 = 1 + (u-1) \, e_i + (u-1) \, e_i \, g_i
\end{equation}
where $g_i$ corresponds to $\sigma_i$ (see \cite{JL} for diagrammatic interpretations). Since the quadratic relations do not change the framing we have $\mathbb{C}C_d^n\subset  {\rm Y}_{d,n}(u)$ and we keep the same notation for the elements of $\mathbb{C}C_d^n$ and for the elements $e_i$  in $ {\rm Y}_{d,n}(u)$. The elements $g_i$ are invertible:
$$
g_i^{-1} = g_i + (u^{-1} - 1)\, e_i +(u^{-1} - 1)\, e_i\, g_i.
$$
For $d=1$ we have $t_j= 1$ and $e_{i} = 1$, and in this case the quadratic relations (\ref{sigmaquadr}) become $g_i^2 = (u-1)g_i +u$, which are the quadratic relations of the Iwahori--Hecke algebra ${\rm H}_n (u)$. So, ${\rm Y}_{1,n}(u)$ coincides with the algebra ${\rm H}_n(u)$.
Further, there is an obvious epimorphism of the Yokonuma-Hecke algebra ${\rm Y}_{d,n}(u)$ onto the algebra ${\rm H}_n(u)$ via the map:
\begin{equation}\label{epimorphism}
\begin{array}{lll} 

g_i & \mapsto & h_i\\
t_j & \mapsto & 1.
\end{array}
\end{equation}

We can alternatively define the algebra ${\rm Y}_{d,n}(u)$ as a $u$--deformation of the algebra $\mathbb{C}C_{d,n}$. More precisely, let $w \in S_n$ and let $w = s_{i_1}\ldots s_{i_k}$ be a reduced expression for $w$. Since the generators $g_i$ of ${\rm Y}_{d,n}(u)$ satisfy the same braiding relations as the generators of $S_n$, then together with the well-known theorem of Matsumoto \cite{ma}, it follows that $g_w: = g_{i_1}\ldots g_{i_k}$ is well defined. Notice that the defining generators $g_i$ correspond to $g_{s_i}$. We have the 
following multiplication rule in ${\rm Y}_{d,n}(u)$ (see Proposition 2.4\cite{jusur}):
\begin{equation}\label{rule}
g_{s_i}g_w =\left\{\begin{array}{ll}
g_{s_iw} & \text{for } l(s_iw)>l(w) \\
 g_{s_iw}+(u-1)e_ig_{s_iw} + (u-1)e_ig_w & \text{for } l(s_iw)<l(w) 
\end{array}\right.
\end{equation}
We also correspond $g_{t_i}$ to $t_i$ and we define: $g_{t_iw} = g_{t_i}g_w =  t_ig_w$. 
Using the above multiplication formulas the second author proved in \cite{Juyu} that ${\rm Y}_{d,n}(u)$ has the following  standard  basis:
$$ 
\{t_1^{a_1}\cdots t_n^{a_n} g_w \, | \,a_i \in {\mathbb Z}/d{\mathbb Z},\,  w \in S_n \}
$$

Further, we have an inductive basis of the Yokonuma--Hecke algebra, which is used in the proof of the main theorem. 

\begin{proposition}[\cite{Juyu} Proposition 8]\label{inductiveYH}
Every element in ${\rm Y}_{d,n+1}(u)$ is a unique linear combination of words, each of one of the following types:
$$
\mathfrak{m}_{n}g_n g_{n-1}\ldots g_i t_i^k \quad \mbox{or} \quad  \mathfrak{m}_{n}t_{n+1}^k, 
$$
where $k\in \mathbb{Z}/d \mathbb{Z}$ and $\mathfrak{m}_n$ is  a word in the inductive basis of  ${\rm Y}_{d,n}(u)$.
\end{proposition}

\subsection{{\it A Markov trace on ${\rm Y}_{d,n}(u)$}}

Using the above basis, the second author constructed in \cite{Juyu} a linear Markov trace on the algebra ${\rm Y}_{d,n}(u)$. Namely:
\begin{theorem}[\cite{Juyu} Theorem 12]\label{Juyutrace}
Let $d$ a positive integer. For indeterminates $z$, $x_1$, $\ldots, x_{d-1}$ there exists a unique linear Markov trace ${\rm tr}$:
$$
 {\rm tr}:  \cup_{n=1}^{\infty}{\rm Y}_{d,n}(u) \longrightarrow   K[z, x_1, \ldots, x_{d-1}]
$$
 defined inductively on $n$ by the following rules:
$$
\begin{array}{rcll}
{\rm tr}(ab) & = & {\rm tr}(ba)  \qquad &  \\
{\rm tr}(1) & = & 1 & \\
{\rm tr}(ag_n) & = & z\, {\rm tr}(a) \qquad &  (\text{Markov  property} )\\
{\rm tr}(at_{n+1}^s) & = & x_s {\rm tr}(a)\qquad  & (  s = 1, \ldots , d-1)
\end{array}
$$
where $a,b \in {\rm Y}_{d,n}(u)$.
\end{theorem}
By direct computation, ${\rm tr}(e_i)$ takes the same value for all $i$. We denote this value by $E$, that is:
$$
E := {\rm tr}(e_i)= \frac{1}{d}\sum_{s=0}^{d-1}x_{s}x_{d-s} ,
$$
where $x_0:=1$. For all $0\leq m \leq d-1$, we also define:

$$
E^{(m)} :={\rm tr}(e_i^{(m)})= \frac{1}{d}\sum_{s=0}^{d-1}x_{m+s}x_{d-s} ,
$$
where $e_i^{(m)}$ is defined in (\ref{eijm}). Notice that  $E = E^{(0)}$.

\subsection{{\it The ${\rm E}$--system.}}

In order for an invariant for framed knots and links to be constructed through the trace on ${\rm Y}_{d,n}(u)$, ${\rm tr}$ should be normalized and rescaled properly. In \cite{JL2} it is proved that such a rescaling is possible if the trace parameters $x_i$ are solutions of a non-linear system of equations, the so--called ${\rm E}$--system.
\begin{definition} \rm
We say that the set of complex numbers $\{x_0,x_1, \ldots ,x_{d-1}\}$ (where $x_0$ is always equal to 1) satisfies the {\it ${\rm E}$--condition} if
$x_1,\ldots , x_{d-1}$ satisfy the following {\it ${\rm E}$--system} of non--linear equations in ${\mathbb C}$:
$$ 
E^{(m)} = x_m E \qquad (1 \leq m \leq d-1)$$
or equivalently:
\begin{equation}\label{Esystem}
\sum_{s=0}^{d-1}x_{m + s} x_{d-s}  =  x_m \sum_{s=0}^{d-1} x_{s} x_{d-s} \qquad (1\leq m \leq d-1).
\end{equation}

\end{definition}

In \cite[Appendix]{JL2} it is proved that the solutions of the  ${\rm E}$--system are the functions $x_{_S}$, from $\mathbb{Z}/d\mathbb{Z}$ to $\mathbb{C}$, parametrized by the non--empty subsets $S$ of the cyclic group $\mathbb{Z}/d\mathbb{Z}$ as follows:
\begin{equation}\label{esystemsol}
x_{_S}={\frac{1}{ |S|}}\sum_{s\in S}\exp_s
\end{equation}
where $ \exp_s(k) = \cos\frac{2\pi s k}{d} + i \sin \frac{2 \pi s k}{d}$ $(k\in \mathbb{Z}/d\mathbb{Z})$. 

\begin{remark} \rm

It is worth noting that the solution of the ${\rm E}$--system can be interpreted as a generalization of the Ramanujan's sum. Indeed, by taking the subset $P$ of  $\mathbb{Z}/d{\mathbb{Z}}$ consisting  
of the numbers coprimes to $d$, then the solution parametrized by $P$ is, up to the factor $\vert P\vert$,  the Ramanujan's sum $c_d(k)$ (see \cite{ra}).
\end{remark}

Equivalently, 
$x_{_S}$ can be seen as an element in $\mathbb{C}C_d$, namely:
\begin{equation}\label{xingrpalg}
x_{_S}= \sum_{k=0}^{d-1} x_k t^k
\end{equation}
where $x_k = \frac{1}{|S|} \sum _{s\in S}\chi_s t^k$, $k=0, \ldots , d-1$, and $\chi_s$ is the character 
of $C_d$ defined as $\chi_s : t^m\mapsto \exp(sm)$. So, the coefficient $x_k$ of $t^k$ in \eqref{xingrpalg} corresponds to $x_{_S}(k)$ in \eqref{esystemsol}.

Recall now that on the group algebra $\mathbb{C}G$ of the finite group $G$,  we have two products, one of them is the multiplication by coordinates, also called the multiplications of the values, which is defined as:
$$
\left( \sum_{g\in G} a_g g\right) \cdot 
\left( \sum_{g\in G} b_g g\right)=
 \sum_{g\in G} a_gb_g g.
$$
and the other product  is the convolution product:
\begin{equation} \label{AA1}
\left( \sum_{g\in G} a_g g\right) * 
\left( \sum_{h\in G} b_h h\right)=
\sum_{g\in G}\sum_{h\in G} a_g b_h gh=
\sum_{g\in G} \left( \sum_{h\in G} a_h b_{gh^{-1}} \right) g.
\end{equation}

By taking $G=C_d$ and writing an arbitrary element  $x$ in $\mathbb{C}C_d$ as $x= \sum_{0\leq k \leq d-1}a_kt^k$,  we have the following lemma:
\begin{lemma}\label{groupalgrels} In $\mathbb{C}C_d$ we have:
$$
x*x= d \sum_{0\leq \ell \leq d-1} E^{(\ell)}t^\ell
$$
and
$$
x*x*x= d^2\sum_{0\leq \ell \leq d-1} \mathrm{tr}(e_1^{\ell}e_2)
 t^\ell.
$$
\end{lemma}

\begin{proof}
The expression for $x\ast x$ follows immediately by direct computation. For the second expression we have that:
\begin{align*}
x\ast x \ast x &= d \sum_{0 \leq \ell \leq d-1} E^{(\ell)} t^\ell \ast x\\
& = d \sum_{0 \leq \ell \leq d-1} E^{(\ell)} t^\ell \ast \sum_{0\leq k \leq d-1}a_kt^k\\
&= d \sum_{0 \leq \ell,k \leq d-1}E^{(\ell)} a_k t^{\ell+k}\\
&= d \sum_{0 \leq \ell,k ,s \leq d-1} a_s a_{\ell -s} a_k t^{\ell +k}\\
&= d \sum_{0 \leq \ell,k ,s \leq d-1} a_s a_{\ell -s-k} a_k t^{\ell}\\
&= d^2 {\rm tr}(e_1^{(\ell)} e_2 ).
\end{align*}
\end{proof}
 For each $a\in \mathbb{Z}/d\mathbb{Z}$ the character $\chi_a$ defines, with respect to the convolution product, an element $\mathbf{i}_a$ of  $\mathbb{C}C_d$,
$$
\mathbf{i}_a:=\sum_{0\leq s\leq d-1} \chi_a(s) t^s.
$$
One can verify that 
$$
\mathbf{i}_a*\mathbf{i}_b=
\begin{cases}
d \, \mathbf{i}_a & \mbox{if } a=b\\
0 & \mbox{if } a\neq b
\end{cases}
$$
that is, $\mathbf{i}_a/d$ is an idempotent element. On the other hand, regarding $\delta_a:= t^a$ as element in $\mathbb{C}C_d$, it is clear 
that, 

$$
\delta_a \cdot \delta_b=
\begin{cases}
\delta_a & \mbox{if } a=b\\
0 & \mbox{if } a\neq b
\end{cases}.
$$
The connection between the two products on $\mathbb{C}C_d$ 
is given by the {\em Fourier transform}. More precisely, the Fourier transform is the linear automorphism on $\mathbb{C}C_d$, defined as: 
\begin{equation}\label{fourtrans}
x:=  \sum_{0\leq r \leq d-1}a_rt^r \mapsto \widehat{x}:= (x\ast \mathbf{i}_s)(0) = \sum_{0\leq \ell \leq d-1} a_\ell\chi_s(d-\ell)
\end{equation}
%where
%$$
%b_s :=\sum_{0\leq \ell \leq d-1} a_\ell\chi_s(d-\ell)
%$$
With the above notation we have:
\begin{lemma}\label{groupalgrels2} The following hold in $\mathbb{C}C_d$:
$$
\widehat{x*y}=\widehat{x}\cdot \widehat{y}, \qquad
\widehat{x\cdot y}= d^{-1}\widehat{x} *\widehat{y},
$$
$$
\widehat{\delta}_a=\mathbf{i}_{-a},\qquad \widehat{\mathbf{i}}_a=d \delta_a,
\qquad  \widehat{\widehat{x\, }}(u)=d x(-u).
$$
\end{lemma}
\begin{proof}
The proof is just a straightforward computation (see \cite{te}).
\end{proof}
\section{The Yokonuma--Temperley--Lieb Algebra}
In this section we define the Temperley--Lieb analogue, in the case of framing, as quotient of ${\rm Y}_{d,n}(u)$ over an appropriate two--sided ideal.
\subsection{\it The Yokonuma--Temperley--Lieb algebra} The Hecke algebra, ${\rm H}_n(u)$,  can be considered as a $u$--deformation of the $\mathbb{C}S_n$, while ${\rm TL}_n(u)$ is the quotient of ${\rm H}_n(u)$ over the two--sided ideal: 

$$
J = \langle h_{i,j} \, ; \, \text{for all}\, i,j \, \text{such that }\, \vert i-j\vert = 1\rangle 
$$
where $h_{i,j}$'s are the Steinberg elements $h_{i,j}:= 1 + h_i + h_{i+1} + h_i h_{i+1} + h_{i+1} h_i + h_i h_{i+1} h_i$. 
It is well-known that that $J$ is a principal ideal. Indeed, 
 $$
 J = \langle h_{1,2} \rangle .
 $$
Notice now that $h_{i,j}$ can be rewritten as 
$$
h_{i,j} = \sum_{\alpha\in W_{i,j}}h_{\alpha}
$$
where $W_{i,j}$ is the subgroup of $S_n$ generated by $s_i$ and $s_j$ (clearly,  $W_{i,j}$ is isomorphic to $S_3$).
On the other hand ${\rm Y}_{d,n}(u)$ can be regarded as a $u$--deformation of $\mathbb{C}[C_d^n \rtimes S_n]$. The symmetric group $S_n$  can be considered as a subgroup of  $C_d^n \rtimes S_n$, therefore  the subgroups $W_{i,j}$ of $S_n$ can be also regarded as subgroups of $C_d^n \rtimes S_n$. Thus, in analogy to the ideal $J$ of  ${\rm H}_n(u)$, it is  natural to consider  the   following ideal $I$ of ${\rm Y}_{d,n}(u)$:
\begin{equation}\label{ideal}
I:= \langle g_{i,j} \, ; \, \text{for all}\, i,j \, \text{such that }\, \vert i-j\vert = 1\rangle 
\end{equation}
where 
\begin{equation}\label{gijdef}
g_{i,j}: = \sum_{\alpha\in W_{i,j}}g_{\alpha} = 1 + g_i + g_{j} + g_i g_{j} + g_{j} g_i + g_i g_{j} g_i .
\end{equation}

We then define:
\begin{definition}\label{ytldef} \rm For $n\geq 3$, the \textit{Yokonuma--Temperley--Lieb} algebra, ${\rm YTL}_{d,n}(u)$, is defined as the quotient:
$$
{\rm YTL}_{d,n}(u)= \frac{{\rm Y}_{d,n}(u)}{ I}.
$$
In other words, the algebra ${\rm YTL}_{d,n}(u)$ can be presented by the generators $ 1, g_1, \ldots, g_{n-1}, t_1,\ldots , t_n$ (by abuse of notation), subject to the following relations:
\begin{align}
g_ig_j &= g_jg_i, \quad |i-j| > 1 \label{YH1}\\
 g_{i+1}g_ig_{i+1}&= g_ig_{i+1}g_i \label{YH2}\\
 g_i^2 &= 1 + (u-1)e_i+ (u-1)e_ig_i \label{quadratic}\\
t_it_j &= t_jt_i,  \quad \mbox{for all } i,j \label{YH4}\\
t_i^{d}  &= 1,\quad  \mbox{for all } i \label{YH3}\\
g_i t_i &=  t_{i+1}g_i \label{YH5}\\
g_it_{i+1} &= t_i g_i \label{YH6}\\
g_it_j &= t_j g_i, \quad \mbox{for } j \neq i, \, \mbox{and } j \neq i+1 \label{YH7}\\
g_ig_{i+1}g_i +g_ig_{i+1}& +g_{i+1}g_i +g_i +g_{i+1} +1 =0\label{YTL}
\end{align}
We shall refer to relations (\ref{YTL}) as \it{the Steinberg  relations}.
\end{definition}
Notice that relations (\ref{YH1})--(\ref{YH7})  are the defining relations of the algebra
${\rm Y}_{d,n}(u)$. Note also that relations (\ref{YTL}) are symmetric with respect to the indices $i$, $i+1$, i.e.:
$$ g_i g_{i+1} g_i = - g_i g_{i+1} - g_{i+1} g_i - g_{i+1} - g_i -1 = g_{i+1} g_i g_{i+1} .$$
 so relations (\ref{YH2}) follow from relations (\ref{YTL}).
\begin{remark}\label{YtoH} \rm
\rm{In analogy to the Yokonuma--Hecke algebra, ${\rm YTL}_{1,n}(u)$ coincides with the algebra ${\rm TL}_n(u)$. Further, the epimorphism (\ref{epimorphism}) induces an epimorphism of the Yokonuma--Temperley--Lieb algebra ${\rm YTL}_{d,n}(u)$ onto the algebra ${\rm TL}_{n}(u)$. Clearly, by relations (\ref{YH5}) and (\ref{YH6}), any monomial in ${\rm YTL}_{d,n}(u)$ inherits the \textit{splitting property} of ${\rm Y}_{d,n}(u)$, that is, it can be written in the form: 
\begin{equation}\label{splittedword}
w= t_{1}^{a_1}\ldots t_{n}^{a_n} g_{i_1}\ldots g_{i_k},
\end{equation}
where: $a_1, \ldots , a_n  \in \mathbb{Z}/ d \mathbb{Z} $.
}
\end{remark}

We shall now prove that $I$ is in fact a principal ideal.
\begin{lemma}\label{lemmaforprop}
The following hold in ${\rm Y}_{d,n}(u)$ for all $i =1, \ldots , n-2$:
$$
\begin{array}{crcl}
(1) & g_i &=& (g_1 \ldots g_{n-1})^{i-1}\, g_1 \,(g_1 \ldots g_{n-1})^{-(i-1)}\\
(2) & g_{i+1}&=& (g_1 \ldots g_{n-1})^{i-1} \, g_2\, (g_1 \ldots g_{n-1})^{-(i-1)}\\
(3) &  g_{i} g_{i+1} &=& (g_1 \ldots g_{n-1})^{i-1} \, g_1g_2\,(g_1 \ldots g_{n-1})^{-(i-1)}\\
(4)& g_{i+1}g_i &=& (g_1 \ldots g_{n-1})^{i-1}\,  g_2 g_1\,(g_1 \ldots g_{n-1})^{-(i-1)}\\
(5)&  g_ig_{i+1} g_i&=& (g_1 \ldots g_{n-1})^{i-1}\,  g_1 g_2 g_1 \,(g_1 \ldots g_{n-1})^{-(i-1)}
\end{array} $$
\end{lemma}

\begin{proof}
We will demonstrate the proof for the cases $(1)$ and $(5)$. The rest of the cases are proved in an analogous manner. For case $(1)$ we have that the statement is true for $i=2$. Indeed:
\begin{align*}
(g_1 \ldots g_{n-1})\,  g_1 \, (g_1\ldots g_{n-1})^{-1} &= g_1 g_2  g_1 g_3 \ldots g_{n-1} (g_1 \ldots g_{n-1})^{-1}\\
&=  g_2 (g_1 g_2\ldots g_{n-1}) (g_1\ldots g_{n-1})^{-1}\\
& = g_2.
\end{align*}
Suppose that the statement is true for $i=k$. We will show that the statement holds for $i=k+1$.
We have:
\begin{align*}
 (g_1 \ldots g_{n-1})^{k} g_1 (g_1 \ldots g_{n-1})^{-k}& =  (g_1 \ldots g_{n-1})(g_1 \ldots g_{n-1})^{k-1}   g_1 (g_1 \ldots g_{n-1})^{-(k-1)}(g_1 \ldots g_{n-1})^{-1}\\
&= (g_1 \ldots g_{n-1}) g_k (g_1 \ldots g_{n-1})^{-1}  \\
&= g_1 \ldots g_{k-1} g_k g_{k+1}  g_k g_{k+2} \ldots g_{n-1} (g_1 \ldots g_{n-1})^{-1} \\
&=  g_1 \ldots g_{k-1} g_{k+1} g_k g_{k+1}  \ldots g_{n-1} (g_1 \ldots g_{n-1})^{-1} \\
&=g_{k+1} ( g_1 \ldots g_{n-1})(g_1 \ldots g_{n-1})^{-1}\\
&= g_{k+1}.
\end{align*} 

For case $(5)$ we have from (1):
\begin{align*}
 g_ig_{i+1}g_i& =(g_1 \ldots g_{n-1})^{i-1} g_1 (g_1 \ldots g_{n-1})^{-(i-1)}(g_1 \ldots g_{n-1})^{i} g_1 (g_1 \ldots g_{n-1})^{-i}\\
&\quad \cdot (g_1 \ldots g_{n-1})^{i-1} g_1 (g_1 \ldots g_{n-1})^{-(i-1)}\\
&= (g_1\ldots g_{n-1})^{i-1}g_1(g_1 \ldots g_{n-1})^{-(i-1)}(g_1\ldots g_{n-1})^{i-1} (g_1\ldots g_{n-1}) \\
&\quad \cdot g_1 (g_1\ldots g_{n-1})^{-1}(g_1\ldots g_{n-1})^{-(i-1)}(g_1 \ldots g_{n-1})^{i-1} g_1 (g_1 \ldots g_{n-1})^{-(i-1)}\\
&= (g_1 \ldots g_{n-1})^{i-1} g_1 g_2 g_1(g_1\ldots g_{n-1})^{-(i-1)}.
\end{align*}

\end{proof}

\begin{corollary}\label{coroldef}
${\rm YTL}_{d,n}(u)$ is the $K$--algebra generated by the set $\{ 1, t_1, \ldots , t_n , g_1 , \ldots , g_{n-1} \}$ whose elements are subject to the defining relations of ${\rm Y}_{d,n}(u)$ and the relation:
$$g_{1,2}= 0. $$
\end{corollary}
\begin{proof}
The result follows using the multiplication rule defined on ${\rm Y}_{d,n}(u)$ and Lemma~\ref{lemmaforprop}.
\end{proof}

\subsection{\it{A presentation with non-invertible generators}}

 In analogy with Eq.~\ref{transf} one can obtain a presentation for the Yokonuma--Temperley--Lieb algebra ${\rm YTL}_{d,n}((u)$ with the non-invertible generators:
\begin{equation}
l_i := \frac{1}{u+1}(g_i+1). \label{ltog}
\end{equation}
In particular we have:
\begin{proposition}\label{noninvert} ${\rm YTL}_{d,n}(u)$ can be viewed as the algebra generated by the elements:
$$1,l_1,\ldots l_{n-1},t_1,\ldots, t_n,$$ which satisfy the following defining relations:
\begin{eqnarray}
t_i^d &=&1, \quad \mbox{for all }i\label{noninveq1}\\
t_it_j&=&t_jt_i, \quad \mbox{for all }i,j\\
l_i t_j &=& t_j l_i, \quad \mbox{for } j\neq i\,\mbox{and } j \neq i+1\\
l_it_i &=& t_{i+1}l_i + \frac{1}{u+1}(t_i - t_{i+1})\\
l_i t_{i+1}&=& t_il_i +\frac{1}{u+1}(t_{i+1}-t_i)\\
l_i^2 &=& \frac{(u-1)e_i+2}{u+1}\,l_i\label{quadrnoninv}\\
l_il_j&=&l_jl_i, \quad |i-j|>1\\
l_il_{i\pm1}l_i &=& \frac{(u-1)e_i +1}{(u+1)^2}\,l_i \label{noninvrel}
\end{eqnarray}
\end{proposition}
\begin{proof}
Obviously, ${\rm YTL}_{d,n}(u)$ is generated by the $l_i$'s and the $t_i$'s. 
It is a straightforward computation to see that relations (\ref{YH1})--(\ref{YTL}) 
are transformed into the relations \eqref{noninveq1} -- (\ref{noninvrel}). However, we shall show here how it works for the quadratic relations \eqref{noninveq1} and the Steinberg relations \eqref{noninvrel}. From Eq.~\ref{ltog} we obtain:
\begin{equation}\label{gtol}
 g_i = (u+1) l_i -1.
\end{equation}
We then have that:
$$g_i^2 = \left ((u+1)^2 l_i -1  \right )^2$$
which is equivalent to:
$$1+ (u-1)e_i + (u-1)e_ig_i = (u+1)^2 l_i^2 - 2 (u-1) l_i +1$$
or equivalently:
$$ (u-1)(u+1)e_i l_i  = (u+1)^2 l_i^2 - 2(u+1)l_i$$
which leads to:
$$ l_i^2 = \frac{(u-1) e_i + 2}{u+1} l_i. $$
which is Eq.~\ref{quadrnoninv}.\\

For the Steinberg elements $g_{i,i\pm1}$ using Eq.~\ref{gtol} we have that:
$$g_{i, i\pm 1}= g_{i}g_{i \pm 1}g_i + g_{i \pm 1}g_i + g_{i} g_{i \pm 1} + g_{i\pm1} + g_i +1 = (u+1)^3 l_i l_{i\pm1} l_i - (u+1)^2 l_i^2 + (u+1) l_i $$
From the Steinberg  relation (\ref{YTL}) and Eq.~\ref{quadrnoninv} we have that:
$$
 (u+1)^2l_i l_{i\pm1}l_i = ((u-1) e_i +1) l_i
$$
or equivalently:
$$
l_il_{i\pm1}l_i = \frac{(u-1)e_i +1}{(u+1)^2}\,l_i,
$$
which is Eq.~\ref{noninvrel}. 
\end{proof}
\begin{remark}{\rm
Setting $d=1$ in the presentation of ${\rm YTL}_{d,n}(u)$ in Proposition~\ref{noninvert}, one obtains the classical presentation of ${\rm TL}_n(u)$, as discussed in Subsection \ref{clasTL}. Note also that, substituting in the braid relation (\ref{YH2}) the $g_i$'s using Eq.~\ref{gtol}, we obtain the equation:
$$ l_i l_{i+1} l_i - \frac{(u-1) e_i +1}{(u+1)^2}\, l_i = l_{i+1} l_i l_{i+1} - \frac{(u-1)e_{i+1} +1}{(u+1)^2} \,l_{i+1} $$
which becomes superfluous, since it can be deduced from Eq.~\ref{noninvrel}. This was to be expected, since the braid relations (\ref{YH2}) were also superfluous.
}
\end{remark}
\section{A spanning set for  The Yokonuma--Temperley--Lieb algebra}
In this section we discuss various properties of a word in ${\rm YTL}_{d,n}(u)$ and we present a spanning set for ${\rm YTL}_{d,n}(u)$ (Proposition~\ref{ytlspan}). Furthermore, using the work of Chlouveraki and Pouchin in \cite{ChPou} we give their formula for the dimension of ${\rm YTL}_{d,n}(u)$ (Proposition~\ref{ytldim}) and we also discuss their results on the linear basis of ${\rm YTL}_{d,n}(u)$ (Theorem~\ref{chpouthm}). We finally compute a basis for ${\rm YTL}_{2,3}(u)$ different than the one of Theorem~\ref{chpouthm}. 
\subsection{{\it}}We have the following definition:
\begin{definition}
In ${\rm YTL}_{d,n}(u)$ we define a length function $l$ as follows:
$$
l(t^a g_{i_1}\ldots g_{i_k}) := l'(s_{i_1}\ldots s_{i_k}),
$$
where $l'$ is the usual {\it length function} of $S_n$ and $t^a : = t_1^{a_1} \ldots t_n^{a_n}\in C_d^n$. A word in ${\rm YTL}_{d,n}(u)$ of the form~(\ref{splittedword}) shall be called \textit{reduced} if it is of minimal length with respect to relations (\ref{YH1})--(\ref{quadratic}), (\ref{YTL}).
\end{definition}
\begin{proposition}\label{LIG}
Each word in ${\rm YTL}_{d,n}(u)$ can be written as a sum of monomials, where the highest and lowest index of the generators $g_i$ appear at most once.
\end{proposition}
\noindent
\begin{proof}  An analogous statement holds for the Yokonuma--Hecke algebra ${\rm Y}_{d,n}(u)$ where only the highest index generators appear at most once \cite[Proposition 8]{Juyu}. Since ${\rm YTL}_{d,n}(u)$ is a quotient of the algebra ${\rm Y}_{d,n}(u)$ the highest index property passes through to the algebra ${\rm YTL}_{d,n}(u)$. The idea is analogous to \cite[Lemma 4.1.2]{JonesIndex} and it is based on induction on the length of reduced words, use of the braid relations and reduction of length using the quadratic relations (\ref{quadratic}). For the case of the lowest index generator $g_i$ we use induction on the length of reduced words and the Steinberg  relations (\ref{YTL}). Indeed, clearly, the statement is true for all words of length $\leq 2$, namely for words of the form $t^a$, $t^a g_1$ $t^ag_1g_2$ and $t^ag_2g_1$.\smallbreak
For words of length 3: Let $w=t^ag_1g_2g_1$. Applying relation (\ref{YH2}) will violate the highest index property of the word, so we must use the Steinberg  relation (\ref{YTL}) and we have:
\begin{eqnarray}
t^ag_1g_2g_1 &=&-t^ag_2g_1 - t^ag_1g_2  -t^ag_2 -t^ag_1-t^a. \nonumber
\label{example}
\end{eqnarray}

\smallbreak
 We assume that the lowest index generator appears at most once in all words of length $\leq$ $r$, and we will show the lowest index property for words of length $r+1$. Let $w=t^a g_{i_1}g_{i_2}\ldots g_{i_k}$ be a reduced word in ${\rm YTL}_{d,n}(u)$ of length $r+1$, and $l ={\rm min} \left \{ i_1 ,\ldots , i_k \right \}$.
\smallbreak

Let first $w=t^aw_1g_lw_2g_lw_3$,  and suppose that $w_2$ does not contain $g_l$. We then have two possibilities:
 \smallbreak
If $w_2$ does not contain $g_{l+1}$, then $g_l$  commutes with all the $g_i$'s in $w_2$ so the length of $w$ can be reduced using the quadratic relations (\ref{quadratic}) for $g_l^2$ and we use the induction hypothesis:
\begin{eqnarray*}
w&=&t^a w_1g_lw_2g_lw_3\\
&=& t^a w_1w_2g_l^2w_3\\
&=&t^a w_1w_2(1+(u-1)e_l+ (u-1)e_lg_l))w_3\\
&=& t^a w_1w_2w_3 + (u-1) t^a w_1w_2e_l w_3 + (u-1)t^a w_1w_2 e_l g_lw_3.
\end{eqnarray*}
If $w_2$ does contain $g_{l+1}$, then, by the induction hypothesis $w_2$ has the form $w_2=v_1 g_{l+1} v_2$, where in $v_1, v_2$ the lowest index generator is at least $g_{l+2}$, hence:
\begin{eqnarray*}
w &=&t^a w_1 g_l v_1 g_{l+1} v_2 g_l w_3\\
&=&t^a w_1v_1g_l g_{l+1} g_l v_2w_3\\
&=&t^a w_1v_1g_{l+1}g_lg_{l+1}v_2w_3,
\end{eqnarray*}
and there is one less occurrence of $g_l$ in $w$. In the case where $l+1=m$, where $m= max \left \{ i_1, \ldots i_k \right \}$,  we apply instead the Steinberg  relation (\ref{YTL}), so no contradiction is caused with respect to the highest index generator.
Continuing in the same manner for all possible pairs of $g_l$ in the word we reduce to having $g_l$ at most once.
\end{proof}
The following proposition gives us a precise spanning set for ${\rm YTL}_{d,n}(u)$.
\begin{proposition}\label{ytlspan}
The following set of reduced words
\begin{eqnarray}
 \Sigma_{d,n} &=& \left\{ t^a(g_{i_1} g_{i_1 - 1} \ldots  g_{i_1 - k_1}) (g_{i_2} g_{i_2-1} \ldots g_{i_2 - k_2})\ldots ( g_{i_p} g_{i_p-1}\ldots g_{i_p - k_p} )\right\}, \label{YTLbasis}
\end{eqnarray}
where $$t^a= t_1^{a_1}\ldots t_{n}^{a_{n}} \in C_d^n, \quad 1 \leq i_1 < i_2 < \ldots <i_p \leq n-1, $$
and $$1 \leq  i_1 - k_1 < i_2 - k_2 < \ldots < i_p - k_p, $$
spans the Yokonuma--Temperley--Lieb algebra ${\rm YTL}_{d,n}(u)$. The highest index generator is $g_{i_p}$ of the rightmost cycle and the lowest index generator is $g_{i_1 - k_1}$ of the leftmost cycle of a word in $\Sigma_{d,n}$.
\end{proposition}
\begin{proof}
 We will prove the statement by induction on the length of a word starting from the linear basis of the Yokonuma--Hecke algebra ${\rm Y}_{d,n}(u)$ \cite[Proposition 8]{Juyu}. Namely,
 \begin{eqnarray}
 \mathcal{B}_{\,{\rm Y}_{d,n}} &=& \left\{ t^{a}(g_{i_1} g_{i_1 -1} \ldots  g_{i_1 - k_1}) (g_{i_2} g_{i_2-1} \ldots g_{i_2 - k_2})\ldots ( g_{i_p} g_{i_p-1}\ldots g_{i_p - k_p} )\right\},\label{YHeckebasis}
\end{eqnarray}
where: $$a \in (\mathbb{Z}/ d \mathbb{Z})^n, \quad 1 \leq i_1 < i_2 < \ldots <i_p \leq n-1, $$ and $ \mathcal{B}_{\,{\rm Y}_{d,n}}$ spans linearly the quotient ${\rm YTL}_{d,n}(u)$ since it is a quotient of ${\rm Y}_{d,n}(u)$. Note that in $ \mathcal{B}_{\,{\rm Y}_{d,n}}$ there is no restriction on the indices $i_1 - k_1, \ldots, i_p-k_p$. Starting now with a word in the set $ \mathcal{B}_{\,{\rm Y}_{d,n}}$, we will show that it is a linear combination of words in the subset $\Sigma_{d,n}$.
The statement holds trivially for words of length 0,1 and 2, since such words are in $\Sigma_{d,n}$. For length 3 consider the representative case of the word $t^a g_1g_2g_1$ which is not in $\Sigma_{d,n}$. Applying the Steinberg  relation (\ref{YTL}) a linear combination of words in $\Sigma_{d,n}$ is obtained (see Eq.~\ref{example}).
Suppose now that the statement holds for all words of length $\leq q$, namely, that any word in $ \mathcal{B}_{\,{\rm Y}_{d,n}}$ of length $q$ can be written as a linear combination of words in $\Sigma_{d,n}$. Let $w$ be a word in $ \mathcal{B}_{\,{\rm Y}_{d,n}}$ of length $q+1$ which is not contained in $\Sigma_{d,n}$. Then $w$ must contain a pair of consecutive cycles:
$$(g_{i_1}g_{i_1 - 1}\ldots g_{k})(g_{i_2}g_{i_2-1}\ldots g_{l}),$$
where $k \geq l$. It suffices to consider the situation where $i_2 =i_1+1$, otherwise the generators of higher index may pass temporarily to the left of the word. Next,  we move the term $g_k$ as far to the right as possible obtaining:
$$(g_{i_1}\ldots g_{k+1})(g_{i_2} \ldots g_{k+2}\underline{g_k g_{k+1} g_k} g_{k-1} \ldots g_l).$$
We now apply the Steinberg  relation (\ref{YTL}) and we obtain five terms, all of length $<q+1$, and we apply the induction hypothesis. More precisely, we have the following five terms: \smallbreak\noindent
$(g_{i_1}\ldots g_{k+1})(g_{i_2}\ldots g_{k+2}\underline{g_{k+1}g_k} g_{k-1}\ldots g_l),$\\
$(g_{i_1}\ldots g_{k+1})(g_{i_2}\ldots g_{k+2}\underline{ g_{k+1}}g_{k-1}\ldots g_l ),$\\
$(g_{i_1}\ldots g_{k+1})(g_{i_2}\ldots g_{k+2}\underline{g_{k}}g_{k-1}\ldots g_l ),$\\
$(g_{i_1}\ldots g_{k+1})(g_{i_2}\ldots g_{k+2} \underline{g_{k}g_{k+1}}g_{k-1}\ldots g_l),$\\
$(g_{i_1}\ldots g_{k+1})(g_{i_2}\ldots g_{k+2}g_{k-1}\ldots g_l ). $
\smallbreak\noindent
To see the exact position of the highest and lowest index generators in the words of $\Sigma_{d,n}$ one can observe that the position of the highest index generator $g_i$ is already clear in the set $ \mathcal{B}_{\,{\rm Y}_{d,n}}$ (cf. \cite{Juyu} \cite{Jones}). To establish the position of the lowest index generator in the words of $\Sigma_{d,n}$ we shall analyze each of the five terms above. In the first term a cycle of smaller length is created and the difference between the lowest indices of the two cycles, $k+1$ and $l$, increases by one, so we need to apply the Steinberg  relation once more and then use the induction hypothesis. In the second term the subword $(g_{k-1} \ldots g_l )$ may pass to the left (since the generator $g_k$ has disappeared), so we obtain the following word:
\begin{equation}\label{leftcycl}
 (g_{k-1} \ldots g_l ) ( g_{i_1} \ldots g_{k+1})(g_{i_2} \ldots g_{k+1}).
\end{equation}
This word contains two cycles with the same lowest index generators, hence we need to apply the Steinberg  relation (\ref{YTL}) and use the induction hypothesis as above. In the third term, $g_k$ returns to its original position and the subword $(g_{k-1} \ldots g_l)$ may pass to the left, obtaining a word in the set $\Sigma_{d,n}$, namely:
\begin{equation}\label{leftcycl2}
(g_{i_1} \ldots g_{k+1}g_k g_{k-1} \ldots g_l)(g_{i_2} \ldots g_{k+2}).
\end{equation}
The same holds for the forth term, which can be rewritten as:
\begin{equation}\label{leftcycl3}
(g_{i_1} \ldots g_{k+1}g_k g_{k-1} \ldots g_l)(g_{i_2} \ldots g_{k+1}).
\end{equation}
Finally, in the fifth term, the subword $(g_{k-1} \ldots g_l)$ may pass to the far left, namely:
\begin{equation}\label{leftcycl4}
(g_{k-1} \ldots g_{l})(g_{i_1} \ldots g_{k+1}) ( g_{i_2}\ldots g_{k+2}),
\end{equation}
which is a word in the set $\Sigma_{d,n}$.
The fact that the lowest index generator $g_i$ appears in the leftmost cycle of the monomial in $\Sigma_{d,n}$ is now clear from (\ref{leftcycl}), (\ref{leftcycl2}), (\ref{leftcycl3}) and (\ref{leftcycl4}).
Concluding, in each application of the Steinberg  relation (\ref{YTL}) the length of $w$ is reduced by at least one, so, from the above and by the induction hypothesis the proof that $\Sigma_{d,n}$ is a spanning set is concluded.
\end{proof}
  M. Chlouveraki and G. Pouchin in \cite{ChPou} have computed the dimension for ${\rm YTL}_{d,n}(u)$ by using the representation theory of the Yokonuma--Hecke algebra \cite{ChLo}. More precisely, they proved the 
  following result.
\begin{proposition}[cf. Proposition 4  \cite{ChPou}]\label{ytldim} The dimension of the Yokonuma--Temperley--Lieb algebra is:
$$ {\rm dim}( {\rm YTL}_{d,n}(u)) = d c_n + \frac{d(d-1)}{2} \sum_{k=1}^{n-1} 	\dbinom{n}{k}^2, $$
where $c_n$ is the $n^{th}$ Catalan number. 
\end{proposition}

\subsection{{\it}} To find an explicit basis for ${\rm YTL}_{d,n}(u)$ Chlouveraki and Pouchin in \cite{ChPou} worked as follows: As mentioned in Remark~\ref{YtoH} each word in ${\rm YTL}_{d,n}(u)$ inherits the splitting property. For each fixed element in the braiding part, they described a set of linear dependence relations among the framing parts (see \cite[Proposition 5]{ChPou}). Using these relations they extracted from $\Sigma_{d,n}$ (recall Eq.~\ref{YTLbasis}) a smaller spanning set for ${\rm YTL}_{d,n}(u)$ and showed that the cardinality of this smaller spanning set is equal to the dimension of the algebra. Thus, it is a basis for ${\rm YTL}_{d,n}(u)$. Before describing this basis, we will need the following notations:

Let $\underline{i}$ and $\underline{k}$ be the following $p$--tuples:
$$ \underline{i} = ( i_1 , \ldots , i_p) \quad \mbox{and} \quad \underline{k} = ( k_1 , \ldots , k_p)$$
and let $\mathcal{I}$ be the set of pairs $(\underline{i}, \underline{k})$ such that:
$$1 \leq i_1 < \ldots < i_p \leq n-1 \quad \mbox{and} \quad 1 \leq i_1 - k_1 < \ldots < i_p - k_p \leq n-1$$
We also denote by $g_{\underline{i}, \underline{k}}$ the element: 
$$g_{\underline{i}, \underline{k}} := (g_{i_1} g_{i_1 - 1} \ldots  g_{i_1 - k_1}) (g_{i_2} g_{i_2-1} \ldots g_{i_2 - k_2})\ldots ( g_{i_p} g_{i_p-1}\ldots g_{i_p - k_p} )$$
Under these notations the set $\Sigma_{d,n}$ can be written as:
$$
\Sigma_{d,n}= \{ t_1^{r_1} \ldots t_n^{r_n}\, g_{\underline{i}, \underline{k}}\, | \, r_1, \ldots , r_n \in \mathbb{Z}/d\mathbb{Z}, (\underline{i} , \underline{k}) \in \mathcal{I} \}. $$ 
The {\it degree of a word} $w = t_1^{r_1} \ldots t_n^{r_n} g_{i_1} \ldots g_{i_m}$ in ${\rm Y}_{d,n}(u)$, denoted $deg(w)$, is defined to be the integer $m$. Set:
$$\Sigma_{d,n}^{<w} : = \{ s \in \Sigma_{d,n} \, | \, deg(s) < deg(w) \}. $$

The group algebra $K(\mathbb{Z}/d\mathbb{Z})^n$ is isomorphic to the subalgebra of ${\rm Y}_{d,n}(u)$ that is generated by the $t_i$'s but not to the subalgebra of ${\rm YTL}_{d,n}(u)$ that is generated by the $t_i$'s. Further, the group algebra $K(\mathbb{Z}/d\mathbb{Z})^n$  has a natural basis, $B_{d,n}$, given by monomials in $t_1, \ldots , t_n$, the following:
$$
B_{d,n} = \{ t_1^{r_1} \ldots t_n^{r_n} \, | \, r_1, \ldots , r_n \in \mathbb{Z}/d\mathbb{Z} \}.
$$
Thus, any element of $K(\mathbb{Z}/d\mathbb{Z})^n$ can be written as a linear combination of words in $B_{d,n}$.
There is a surjective algebra morphism from $K(\mathbb{Z}/d\mathbb{Z})^n$  to the subalgebra of ${\rm YTL}_{d,n}(u)$ that is generated by the $t_i$'s. We will denote the image of an element $b\in B_{d,n}$ into the subalgebra of ${\rm YTL}_{d,n}(u)$ that is generated by the $t_i$'s with $\bar{b}$. We then have the following theorem:
\begin{theorem}[Chlouveraki and Pouchin, cf. \cite{ChPou}, Theorem~2]\label{chpouthm}
The following set is a linear basis for ${\rm YTL}_{d,n}(u)$:
$$ S_{d,n} = \{ \bar{b}_{\underline{i}, \underline{k}} \, g_{\underline{i}, \underline{k}} \, | \, (\underline{i}, \underline{k}) \in \mathcal{I}, \, b_{\underline{i}, \underline{k}} \in \mathcal{B}_{d,n}(g_{\underline{i}, \underline{k}})\},
$$
where $\mathcal{B}_{d,n}(g_{\underline{i}, \underline{k}})$ is a proper subset of $B_{d,n}$ such that:
$$\{  b_{\underline{i}, \underline{k}} + R(g_{\underline{i}, \underline{k}}) \, | \, b_{\underline{i}, \underline{k}} \in \mathcal{B}(g_{\underline{i}, \underline{k}})\}
$$
is a basis of the quotient space $K(\mathbb{Z}/d\mathbb{Z})^n / R(g_{\underline{i}, \underline{k}})$, and where $R(w)$ is the following ideal of $K(\mathbb{Z}/d\mathbb{Z})^n$:
$$R(w) = \{ \mathfrak{m} \in K(\mathbb{Z}/d\mathbb{Z})^n \, | \, \bar{\mathfrak{m}}\,w \in Span_{\mathbb{C}(u)}(\Sigma_{d,n}^{<w})\}. $$
\end{theorem}
\subsection{{\it}}For $d=2$, $n=3$ it is relatively easy to find a basis for ${\rm YTL}_{2,3}(u)$. We will give here a basis different than the one in Theorem~\ref{chpouthm}.  
Before continuing, we need the following technical  lemma that will be also used in the proof of Theorem~\ref{generalcase}. 
\begin{lemma}[cf. Lemma 7.5 \cite{juptl}]\label{condeqs}
For the element $g_{1,2}$ we have in ${\rm Y}_{d,n}(u)$ (recall (\ref{eij}) for $e_{1,3}$):
$$
\begin{array}{crcl}
(1) &g_1 g_{1,2} &=& [1+ (u-1) e_1]g_{1,2}\\
(2) &g_2g_{1,2} &=& [1 +(u-1)e_2]g_{1,2}\\
(3) &g_1 g_2 g_{1,2} &=& [ 1 +(u-1) e_1+(u-1) e_{1,3}+(u-1)^2 e_1e_2 ]g_{1,2}\\
(4) &g_2g_1g_{1,2}&=& [ 1 +(u-1)\, e_2+(u-1) e_{1,3}+(u-1)^2e_1e_2 ] g_{1,2}\\
(5) &g_1g_2g_1g_{1,2}&=& [ 1+ (u-1)(e_1+e_2+e_{1,3})+(u-1)^2(u+2)\,e_1e_2 ] g_{1,2}
\end{array}
$$
Analogous relations hold for multiplications with $g_{1,2}$ from the right.
\end{lemma}
\begin{proof}
The idea is to expand the left--hand side of each equation and then use Eq.~\ref{quadratic} and Lemma~\ref{eipropop}. We will demonstrate the proof for the indicative cases $(1)$ and $(4)$. The other cases are proved similarly.
\smallbreak
For case $(1)$ we have:
\begin{eqnarray*}
g_1g_{1,2} &=& g_1 + g_1^2 + g_1g_2 + g_1^2 g_2 +g_1g_2g_1 + g_1^2g_2g_1\\
&=& g_1 + \left[ 1 + (u-1) e_1 + (u-1)e_1 g_1  \right]  \\
&& +g_1g_2 +\left[ g_2 +(u-1)e_1g_2 + (u-1)e_1g_1g_2 \right] \\
&&+ g_1g_2g_1 + \left[ g_2g_1  + (u-1)e_1g_2g_1 + (u-1)e_1g_1g_2g_1 \right] \\
&=& g_{1,2} + (u-1)e_1g_{1,2}.
\end{eqnarray*}
Case $(2)$ is completely analogous. In order to prove case $(4)$ we will use cases $(1)$ and $(2)$:
\begin{eqnarray*}
g_2g_1g_{1,2} &=& g_2 \left( g_{1,2} + (u- 1)e_1 g_{1,2}\right)\\
&=& g_2 g_{1,2} + (u-1)e_{1,3}g_2 g_{1,2} \quad \mbox{(Lemma~\ref{eipropop})}\\
&=& \left[ 1+(u-1)e_2 \right] g_{1,2} +(u-1)e_{1,3}(1+(u-1)e_2)g_{1,2}\\
&=&\left[ 1+(u-1)e_2\right] g_{1,2} + (u-1)e_{1,3}g_{1,2} + (u-1)^2e_{1,3}e_2g_{1,2} \quad \mbox{(Lemma~\ref{eipropop})}\\
&=&\left[ 1 +(u-1)\, e_2+(u-1) e_{1,3}+(u-1)^2e_1e_2\right] g_{1,2}.
\end{eqnarray*}

\end{proof}
To find a basis for ${\rm YTL}_{2,3}(u)$: From Proposition~\ref{ytldim} we have that ${\rm dim}({\rm YTL}_{2,3}(u)) = 28$. On the other hand the spanning set $\Sigma_{2,3}$ of ${\rm YTL}_{2,3}(u)$ of Proposition~\ref{ytlspan}, contains 40 elements. Thus, any relation $w_1 g_{1,2} w_2=0$ with $w_1,\, w_2 \in {\rm Y}_{2,3}(u)$ reduces to having $w_1,\,w_2 \in \Sigma_{2,3}$. Further, if any of $w_1, \, w_2$ contain braiding generators, then by Lemma~\ref{condeqs} (after pushing framing generators in $w_2$ to the right) these get absorbed by $g_{1,2}$. Thus, and since $e_{i,j}= \frac{1}{2} (1 +t_i t_{j})$ for $d=2$, it suffices to consider the following system of equations:
\begin{equation}\label{basiseqs}
w_1 \, g_{1,2} w_2 = 0 \qquad  w_1, w_2 \in \mathcal{T},
\end{equation}
where $\mathcal{T} := \{ 1 , t_1 , t_2 , t_3 , t_1 t_2 ,t_1 t_3 ,$ $t_2 t_3, t_1 t_2 t_3 \} $. For finding all possible linear dependencies in $\Sigma_{2,3}$, after substituting  $g_1 g_2 g_1$ with  $-1 -g_1 - g_2 -g_1 g_2 - g_2 g_1$ in Eq.~\ref{basiseqs}, note that some of these 64 equations reduce trivially to $g_{1,2}=0$; for example if $w_2=1$ or $w_2 = t_1 t_2 t_3$ (since it commutes with $g_{1,2}$). From the rest one can extract 12 linearly independent equations which, applied on the spanning set $\Sigma_{2,3}$ lead to the following basis for ${\rm YTL}_{2,3}(u)$:
\begin{eqnarray*}
&&\mathcal{S}_{2,3} = \left \{ 1, t_1 , t_2 , t_1 t_2, g_1, t_2 g_1 , t_3 g_1, t_2 t_3 g_1, g_2 , t_1 g_2 , t_3 g_2 , t_1 t_3 g_2 ,  \right.\\ 
&& g_1 g_2, t_1 g_1g_2, t_2 g_1 g_2 , t_3 g_1 g_2 ,t_1 t_2 g_1 g_2 , t_1 t_3 g_1 g_2 , t_2 t_3 g_1 g_2 , t_1 t_2 t_3 g_1 g_2, \\
 &&\left.  g_2 g_1 , t_1 g_2 g_1 , t_2 g_2 g_1, t_3 g_2 g_1, t_1 t_2 g_2 g_1 , t_1 t_3 g_2 g_1 , t_2 t_3 g_2 g_1, t_1 t_2 t_3 g_2 g_1 \right \}.
 \end{eqnarray*}
 
 \section{A Markov trace on ${\rm YTL}_{d,n}(u)$}\label{ytlsection}
The following section is dedicated to finding the necessary and sufficient conditions for the trace ${\rm tr}$ on ${\rm Y}_{d,n}(u)$ to pass to the quotient algebra ${\rm YTL}_{d,n}(u)$, in analogy to the classical case, where the Ocneanu trace on ${\rm H}_n(u)$ passes to the quotient algebra ${\rm TL}_n(u)$ under the condition for certain values of the trace parameter $\zeta$. 
\smallbreak

\subsection{{\it}}It is clear by now that ${\rm tr}$ will pass to ${\rm YTL}_{d,n}(u)$ if it kills the generator of the principal ideal through which the quotient is defined, that is, if ${\rm tr}(g_{1,2})=0$. We have the following lemma:
\begin{lemma}\label{trg12}
For the element $g_{1,2}$ we have:
\begin{equation}\label{trg12eq}
{\rm tr}(g_{1,2}) = (u+1)z^2 + \left((u-1)E +3 \right)z +1.
\end{equation}
\end{lemma}
\begin{proof}
The proof is a straightforward computation:
\begin{align*}
{\rm tr}(g_{1,2})& = {\rm tr}(1) + {\rm tr}(g_1) + {\rm tr}(g_2) +{\rm tr}(g_1g_2) + {\rm tr}(g_2g_1) + {\rm tr}(g_1g_2g_1)\\
&= 1 + 2z + 2z^2 + z+  (u-1) E z  + (u-1)z^2\\
& = (u+1)z^2 + \left ((u-1)E + 3\right)z + 1.
\end{align*}
\end{proof}
Lemma~\ref{trg12}, together with the equation: 
\begin{equation}\label{trg120}
{\rm tr}(g_{1,2})= (u+1)z^2 + \left ((u-1)E + 3\right)z + 1 =0
\end{equation}
 gives us the following values for $z$:
\begin{eqnarray}\label{zval}
z_{\pm}= \frac{ - \left((u -1) E + 3 \right) \pm \sqrt{\left((u-1) E +3\right)^2 - 4(u+1)} }{2(u+1)}.
\end{eqnarray}

We shall do now  the analysis for all conditions that must be imposed on the trace parameters in order that  $\rm tr$ passes to ${\rm YTL}_{d,n}(u)$. Having in mind Corollary \ref{coroldef} and the linearity of ${\rm tr}$, it follows that $\rm tr$ passes to ${\rm YTL}_{d,n}(u)$ if and only if  the following equations are satisfied for all monomials $\mathfrak{m}$ in the inductive basis of ${\rm Y}_{d,n}(u)$. Namely:
\begin{equation}\label{trwg12}
{\rm tr}(\mathfrak{m}\, g_{1,2})=0 .
\end{equation}

Let us first consider the case $n=3$. By Proposition~\ref{inductiveYH} the elements in the inductive basis of   ${\rm Y}_{d,3}(u)$ are of the following forms:
\begin{equation}\label{basicwords}
t_1^{a}t_2^{b}t_3^c,  \quad t_1^{a}g_1t_1^{b}t_3^c, \quad t_1^{a}t_2^{b}g_2g_1t_1^c, \quad  t_1^{a}t_2^{b}g_2t_2^c,\quad t_1^{a}g_1t_1^{b} g_2t_2^c, \quad t_1^{a}g_1t_1^{b} g_2g_1t_1^c
\end{equation}

Using Lemma~\ref{condeqs} and the following notations:
\begin{align*}
& Z_{a,b,c}: = (u+1)z^2 x_{a+b+c} + \left( (u-1) E^{(a+b+c)} + x_{a}x_{b+c} + x_{b} x_{a+c} + x_{c} x_{a+b} \right)z + x_{a} x_{b} x_{c}\\
&V_{a,b+c}:=(u+1)z^2 x_{a+b+c} + (u+1)zE^{(a+b+c)} + z\,x_{a}x_{b+c} + x_{a}E^{(b+c)} \\
&V_{b,a+c}:=(u+1)z^2 x_{a+b+c} + (u+1)zE^{(a+b+c)} + z\,x_{b}x_{a+c} + x_{b}E^{(a+c)} \\
&V_{c,a+b}:=(u+1)z^2 x_{a+b+c} + (u+1)zE^{(a+b+c)} + z\,x_{c}x_{a+b} + x_{c}E^{(a+b)} \\
&W_{a,b,c}:= (u+1)z^2 x_{a+b+c} + (u+2)zE^{(a+b+c)} + {\rm tr}\left(e_1^{(a+b+c)}e_2 \right)  
\end{align*}
we obtain by \eqref{trwg12} and \eqref{basicwords}  the following equations, for any $a,b,c \in \mathbb{Z}/d\mathbb{Z}$:
\begin{align}
& Z_{a,b,c} = 0\label{ytleq1} \\
&Z_{a,b,c}+ (u-1) V_{c,a+b} =0 \label{ytleq2} \\
&Z_{a,b,c} + (u-1)V_{a,b+c} = 0 \label{ytleq3}\\
& Z_{a,b,c}+ (u-1) \left [V_{c,a+b} + V_{b,a+c} + W_{a,b,c} \right ]=0 \label{ytleq4}\\
& Z_{a,b,c}+ (u-1) \left [V_{a,b+c} + V_{b,a+c} + W_{a,b,c} \right ]=0 \label{ytleq5}\\
& Z_{a,b,c} + (u-1) \left [V_{a,b+c} + V_{b,a+c} + V_{c,a+b}+ W_{a,b,c} \right ]=0 \label{ytleq6}
\end{align}
Equations~\ref{ytleq1}--\ref{ytleq6} reduce to the following system of equations of $z,x_1, \ldots , x_{d-1}$ for any $a,b,c \in \mathbb{Z}/d\mathbb{Z}$:
\begin{subnumcases}{(\Sigma)}
Z_{a,b,c} = 0\label{redeq1}\\
V_{c,a+b} =0\label{redeq2}\\
 V_{a,b+c} = 0\label{redeq3}\\
 V_{b,a+c} + W_{a,b,c} =0\label{redeq4} 
\end{subnumcases}
Notice that for $a=b=c=0$ Eq.~\ref{ytleq1} becomes Eq.~\ref{trg120}. If, now, we require both solutions in \eqref{zval} to participate in the solutions of $(\Sigma)$, then we are led to necessary conditions for ${\rm tr}$ to pass to ${\rm YTL}_{2,3}(u)$ (Section 4.2). If not then we are led to necessary and sufficient conditions for ${\rm tr}$ to pass to ${\rm YTL}_{2,3}(u)$ (Section 4.3).

\subsection{{\it}}

Suppose that both solutions for $z$ from Eq.~\ref{zval} participate in the solution set of ($\Sigma$). We have the following proposition:
\begin{proposition}\label{ytlreqtrace}
The trace ${\rm tr}$ defined on ${\rm Y}_{d,3}(u)$ passes to the quotient ${\rm YTL}_{d,3}(u)$ if the trace parameters $x_i$ are $d^{th}$ roots of unity ($x_i= x_1^i$, $ 1 \leq i \leq d-1$) and $z= -\frac{1}{u+1}$ or $z=-1$.
\end{proposition}
\begin{proof}
Suppose that ${\rm tr}$ passes to ${\rm YTL}_{d,3}(u)$ and that $(\Sigma)$ has both solutions for $z$ from Eq.~\ref{zval}. This implies that there exist  $\lambda$ in $K( x_1 \ldots, x_{d-1})$ such that:
$$Z_{a,b,c} = \lambda Z_{0,0,0} $$
From this we deduce that: 
\begin{align}
\lambda &= x_{a+b+c}\nonumber\\
x_ax_{b+c} + x_b x_{a+c} + x_c x_{a+b} &= 3 x_{a+b+c}\nonumber\\
E^{(a+b+c)}& = x_{a+b+c} E\label{ytltrcond1}\\
x_{a+b+c}&=x_a x_b x_c\label{ytltrcond2}.
\end{align}
Since this holds for any $a,b,c \in \mathbb{Z}/d\mathbb{Z}$, by taking $b=c=0$ in Eq.~\ref{ytltrcond1} we have that:
\begin{equation}\label{Eabc}
E^{(a)} = x_a E
\end{equation}
which is exactly the ${\rm E}$--system. Moreover, by taking $c=0$ in Eq.~\ref{ytltrcond2} we obtain:
\begin{equation}\label{rootsofunity}
x_a x_b = x_{a+b}
\end{equation}
This implies that the $x_i$'s are $d^{th}$ roots of unity which is equivalent to $E=1$ \cite[Appendix]{JL2}.  
In order to conclude the proof it is enough to verify that these conditions for the $x_i$'s satisfy also \eqref{redeq2}--\eqref{redeq4} of $(\Sigma)$. Since the $x_i$'s are solutions of the ${\rm E}$--system, Eqs.~\ref{redeq2} and \ref{redeq3} are immediately satisfied. We will finally check Eq.~\ref{redeq4}.
Using Eqs.~\ref{rootsofunity} and \ref{Eabc} we have that:
$$u \left ( (u+1)z^2 + (u+2)z + 1 \right ) x_a x_b x_c = 0,$$
from which we deduce that $ z= - \frac{1}{u+1}$ or $z=-1$, which are precisely the solutions \eqref{zval} for $E=1$.
\end{proof}
Using induction on $n$ one can prove the general case of the necessary conditions for ${\rm tr}$ to pass to ${\rm YTL}_{d,n}(u)$. Indeed we have:
\begin{theorem}\label{generalcase}
For $n \geq 3$, the trace ${\rm tr}$ defined on ${\rm Y}_{d,n}(u)$ passes to the quotient ${\rm YTL}_{d,n}(u)$ if the trace parameters $x_i$ are $d^{th}$ roots of unity ($x_i= x_1^i$, $ 1 \leq i \leq d-1$) and $z= -\frac{1}{u+1}$ or $z=-1$.
\end{theorem}
\begin{proof}
By induction on $n$. In Proposition~\ref{ytlreqtrace} we proved the case where $n=3$. Assume that the statement holds for all ${\rm YTL}_{d,k}(u)$, where $ k \leq n$, that is:
$$
{\rm tr} ( a_k \, g_{1,2} ) = 0
$$
for all $a_k \in {\rm Y}_{d,k}(u)$, $k \leq n$. We will show the statement for $k=n+1$. It suffices to prove that the trace vanishes on any element in the form
 $a_{n+1}g_{1,2}$, where $a_{n+1}$ belongs to the inductive basis of ${\rm Y}_{d,n+1}(u)$ (recall Proposition~\ref{inductiveYH}), given the conditions of the Theorem. Namely:
$$ {\rm tr}(a_{n+1}\, g_{1,2}) =0.$$
Since $a_{n+1}$ is in the inductive basis of ${\rm Y}_{d,n+1}(u)$, it is of one of the following forms:
$$
a_{n+1}=a_{n} g_{n} \ldots g_i t_i^k \quad \mbox{or} \quad a_{n+1}=a_{n}t_{n+1}^k,
$$
where $a_n$ is in the inductive basis of ${\rm Y}_{d,n}(u)$. For the first case we have:
$$
{\rm tr} (a_{n+1} \, g_{1,2}) = {\rm tr} (a_{n} g_{n} \ldots g_i t_i^k \, g_{1,2})
= z\, {\rm tr}(a_{n} g_{n-1} \ldots g_i t_i^k\, r_{1,2}) =z\, {\rm tr}( \tilde{a} g_{1,2}),
$$
where $\tilde{a}:=a_{n} g_{n-1} \ldots g_i t_i^k$. Notice now that $\widetilde{a} $ is a word in ${\rm Y}_{d,n}(u)$ and so, by the linearity of the trace, we have that ${\rm tr}(\tilde{a} \,g_{1,2})$ is a linear combination of traces of the form ${\rm tr}(a_{n}\,g_{1,2})$, where $a_{n}$ is in the inductive basis of ${\rm Y}_{d,n}(u)$. Therefore, by the induction hypothesis, we deduce that:
$$ 
{\rm tr}(\widetilde{a} \, g_{1,2}) =0,
$$
 if the conditions of the Theorem are satisfied. Therefore the statement is proved. The second case is proved similarly. Hence, the proof is concluded.
\end{proof}
\subsection{{\it}}
In the proofs of Proposition~\ref{ytlreqtrace} and Theorem~\ref{generalcase} it became apparent that the $x_i$'s are $d^{th}$ roots of unity if and only if the values of $z_+$ and $z_{-}$ satisfy all equations of $(\Sigma)$. Clearly, if we loosen this last condition, then other solutions for the $x_i$'s may appear such that the trace ${\rm tr}$ passes to the quotient ${\rm YTL}_{d,n}(u)$.  Indeed, we have the following:
\begin{theorem}\label{ytlthmgen}
The trace ${\rm tr}$ passes to the quotient ${\rm YTL}_{d,n}(u)$ if and only if the $x_i$'s  are solutions of the ${\rm E}$--system and one of the two cases holds: 
\begin{enumerate}
\item [(i)] 
For some  $0 \leq m_1 \leq d-1$ the $x_\ell$'s are expressed as:
\[
x_\ell=\exp(\ell m _1) \quad (0 \leq \ell \leq d-1).
\]
In this case  the $x_\ell$'s are $d^{th}$ roots of unity and 
$z=-\frac{1}{u+1}$ or $z=-1$.
\item [(ii)] 
For some $0 \leq m_1,m_2 \leq d-1$ the $x_\ell$'s are expressed as: 
$$
x_\ell=\frac{1}{2}\left(\exp(\ell m _1) +\exp(\ell m _2) \right) \quad (0 \leq \ell \leq d-1).
$$
In this case 
 we have $z=-\frac{1}{2}$.
\end{enumerate}
%
%
%$$x=  \sum_{m \in S_1} \mathbf{i}_{-m } +
%\frac{1}{2}\sum_{m \in S_{\frac{1}{2}}} \mathbf{i}_{-m}
%$$
%under the condition: $$x(0)=1$$ where $S_1$, $S_{\frac{1}{2}}$ are the subsets of the partition of the set $\left \{ \ell \, : \, 0 \leq \ell \leq d-1 \right \}$ into three subsets, $S_0$, $S_1$ and $S_{\frac{1}{2}}$, which correspond to the solutions $y_\ell =0$, $y_\ell = d$ and $y_\ell = \frac{1}{2} d$ of the equation:
%$$y_\ell \left (\frac{2}{d^2} y_\ell^2 -\frac{3}{d}y_\ell+1 \right)=0.$$
%
%More precisely:
\end{theorem}
Note that case (i) captures Theorem~\ref{generalcase}.
\begin{proof}
 Observe that the $x_\ell$'s expressed by (i) are indeed solutions of the system $(\Sigma)$. We will now assume that our solutions are not of this form. This implies that $x_\ell\neq E^{(\ell)}$
for some $0\leq \ell \leq d-1$, and this will allow us to have this quantity in denominators later.

We will use  induction on $n$. We will first prove the case $n=3$. Suppose that trace ${\rm tr}$ passes to the quotient algebra ${\rm YTL}_{d,3}(u)$. This means that $(\Sigma)$ has solutions for $z$ any one of those in Eq.~\ref{zval}, for any $a,b,c \in \mathbb{Z}/d\mathbb{Z}$. Subtracting Eq.~\ref{redeq1} from Eq.~\ref{redeq2} we obtain: 
\begin{equation}\label{zval1/2}
z= - \frac{x_ax_bx_c - x_c E^{(a+b)}}{x_ax_{b+c} + x_bx_{a+c} - 2 E^{(a+b+c)}}.
\end{equation}

For $b=c=0$ in Eq.~\ref{zval1/2} we obtain: $z = - \frac{1}{2}$.  On the other hand, subtracting Eqs.~\ref{redeq1} and \ref{redeq2} from Eq.~\ref{redeq4} we have:
\begin{equation}\label{zvalabc}
z= \frac{x_ax_bx_c + x_c E^{(a+b)} - x_bE^{(a+c)} - {\rm tr}(e_1^{(a+b+c)}e_2)}{3E^{(a+b+c)} - x_ax_{b+c} -2 x_cx_{a+b}}.
\end{equation}
We will now assume that that $x_\ell$'s are not roots of unity. This implies that, for all $0\leq a \leq d-1$, $x_a-E^{(a)}\neq 0$. For $b=c=0$ in Eq.~\ref{zvalabc} we obtain:
\begin{equation}\label{zval1/3x}
z= - \frac{x_a - {\rm tr}(e_1^{(a)} e_2)}{3 ( x_a -  E^{(a)})}.
\end{equation}
From Eqs.~\ref{zval1/2} and \ref{zval1/3x} we have that:
$$
\frac{1}{2} =  \frac{x_a - {\rm tr}(e_1^{(a)} e_2)}{3 ( x_a -  E^{(a)})}
$$
or equivalently:
$$3(x_a - E^{(a)}) = 2( x_a - {\rm tr}(e_1^{(a)}e_2)). $$
Using Lemma~\ref{groupalgrels}, this is equivalent to: 
$$
3x - \frac{3}{d} x*x = 2x - \frac{2}{d^2} x*x*x.
$$
By taking the Fourier transform (see Lemma~\ref{groupalgrels2}) we arrive at:
$$
\frac{2}{d^2}\widehat{x}^3-\frac{3}{d}\widehat{x}^2+\widehat{x}=0.
$$
Assuming that $\widehat{x}=\sum_{0 \leq  \ell \leq d-1 } y_\ell t^\ell$ we have the following expression for the cofficients $y_\ell$ in the expansion of $\widehat{x}$:
$$
y_\ell \left (\frac{2}{d^2} y_\ell^2 -\frac{3}{d}y_\ell+1\right)=0.
$$
So either $y_\ell=0$ or $y_\ell=d$ or $y_\ell=\frac{1}{2}d$. 
So if we take a partition of the set $\{\ell: 0\leq \ell\leq d-1\}
$ into sets $S_0$, $S_1$, $S_{\frac{1}{2}}$ such that $y_\ell$
takes the value $i\cdot d$ on $S_i$ $(i=0,1,\frac{1}{2})$. 
We have from Lemma~\ref{groupalgrels2} that: 
$$
x=  \sum_{m  \in S_1} \mathbf{i}_{-m} +
\frac{1}{2}\sum_{m  \in S_{\frac{1}{2}}} \mathbf{i}_{-m}.
$$
From $x_0=1$ we obtain the conditions:
$$
1=x(0)= |S_1| + \frac{1}{2} |S_{\frac{1}{2}}|.
$$

This means that either $S_1$ has only one element and $S_{\frac{1}{2}}=\emptyset$
or $S_1=\emptyset$ and $S_{\frac{1}{2}}$ has two elements. 
The first case corresponds to the case ${\rm (i)}$ where the $x_\ell$'s are $d^{th}$ 
roots of unity. In the second case, if $S_{\frac{1}{2}}=\{m_1,m_2\}$ we obtain the following solution of the ${\rm E}$--system:

$$
x_\ell=\frac{1}{2}\left(\exp(\ell m _1) +\exp(\ell m _2) \right), \quad (0 \leq \ell \leq d-1)
$$
which corresponds to $z= - \frac{1}{2}$.\\
 The rest of proof (the induction on $n$) is analogous to the one of Theorem~\ref{generalcase}.
\end{proof}
\begin{remark} \rm
The values for the trace parameter $z$ in Theorems~\ref{generalcase} and \ref{ytlthmgen}, $z=-\frac{1}{u+1}$ and $z=-1$, in order that ${\rm tr}$ on ${\rm Y}_{d,n}(u)$ passes to the quotient ${\rm YTL}_{d,n}(u)$ are the same as the values in Eq.~\ref{jonval} for $\zeta$ of the Ocneanu trace $\tau$ on ${\rm H}_n(u)$, so that $\tau$ passes to the quotient ${\rm TL}_n(u)$ (recall Section 1.2).
\end{remark}
\section{Knot Invariants from ${\rm YTL}_{d,n}(u)$}
The 2--variable Jones or HOMFLYPT polynomial, $P(\lambda, u)$, can be defined through the Ocneanu trace on ${\rm H}_n(u)$ \cite{Jones}. Indeed, for any braid $\alpha \, \in \, \cup_{\infty} B_{n}$ we have:
$$ P(\lambda,u)(\hat{\alpha}) = \left( - \frac{1-\lambda u}{\sqrt{\lambda}(1-u)}\right)^{n-1} \left (\sqrt{\lambda}\right)^{\varepsilon(\alpha)} {\rm \tau}(\pi(\alpha)),$$
where: $\lambda = \frac{1-u+\zeta}{u\zeta}$, $\pi$ is the natural epimorphism of $\mathbb{C}B_n$ onto ${\rm H}_n(u)$ that sends the braid generator $\sigma_i$ to $h_i$ and $\varepsilon(\alpha)$ is the algebraic sum of the exponents of the $\sigma_i$'s in $\alpha$. Further, the Jones polynomial, $V(u)$, related to the algebras ${\rm TL}_n(u)$, can be redefined through the HOMFLYPT polynomial, related the algebras ${\rm H}_n(u)$, by specializing $\zeta$ to $-\frac{1}{u+1}$ \cite{Jones}. This is the non--trivial value for which the Ocneanu trace $\tau$ passes to the quotient ${\rm TL}_n(u)$. Namely:
$$V(u)(\hat{\alpha}) = \left(- \frac{1+u}{\sqrt{u}} \right)^{n-1} \left(\sqrt{u}\right)^{\varepsilon(\alpha)} {\rm \tau}(\pi(\alpha)) = P(u,u)(\hat{\alpha}).$$

\subsection{{\it}} In \cite{JL2} it is proved that the trace ${\rm tr}$ can be re-scaled according to the braid equivalence corresponding to isotopic framed links if and only if the $x_i$'s furnish a solution of the ${\rm E}$--system. Then, by further normalizing an invariant for framed knots and links can be obtained \cite{JL2}:
\begin{equation}\label{gammainv}
\Gamma_{d,S}(w,u)(\hat{\alpha}) = \left(- \frac{1 - w u}{\sqrt{w} (1-u) E} \right)^{n-1} \left(\sqrt{w}\right)^{\varepsilon(\alpha)} {\rm tr}(\gamma(\alpha)),
\end{equation}
where: $S$ is a subset of $\mathbb{Z}/d\mathbb{Z}$ which parametrizes a solution of the ${\rm E}$--system, $w = \frac{z + (1-u)E}{uz}$, $\gamma$ the natural epimorphism of the framed braid group algebra $\mathbb{C}\mathcal{F}_n$ onto the algebra ${\rm Y}_{d,n}(u)$, and $\alpha \in \cup_{\infty} \mathcal{F}_{n}$. Note that  for every $d \in \mathbb{N}$ we obtain $2^d -1 $ invariants for framed links.

 Further, in \cite{JL4}  the second and the fourth authors represented the classical braid group $B_n$ in the algebra ${\rm Y}_{d,n}(u)$ by regarding the framing generators $t_i$ as formal elements. So, $\Gamma_{d,S}$ can be seen as an invariant of classical links. Namely:
\begin{equation}\label{deltainv}
\Delta_{d,S}(w,u)(\hat{\alpha}) = \left(- \frac{1 - w u}{\sqrt{w} (1-u) E} \right)^{n-1} \left(\sqrt{w}\right)^{\varepsilon(\alpha)} {\rm tr}(\delta(\alpha)),
\end{equation}
where: $S, \,w$ as above, $\delta$ the natural homomorphism of the classical braid group algebra $\mathbb{C}B_n$ to the algebra ${\rm Y}_{d,n}(u)$ and $\alpha \in \cup_{\infty}B_n$. Further, in \cite{JL3} the invariant $\Delta_{d,S}(w,u)$ was extended to an invariant for singular links. 

In \cite{ChLa} it is shown that for generic values of the parameters $u,z$ the invariants $\Delta_{d,S}(w,u)$ do not coincide with the HOMFLYPT polynomial except in the trivial cases $u=1$ or $E=1$. Yet, computational data \cite{CJJKL} indicate that these invariants do not distinguish more or less knot pairs than the HOMFLYPT polynomial, so they may still be topologically equivalent to the HOMFLYPT polynomial. 

\subsection{{\it}}We shall now define framed and classical link invariants related to the algebra ${\rm YTL}_{d,n}(u)$. In Theorem~\ref{ytlthmgen} we showed that the trace ${\rm tr}$ passes to the quotient ${\rm YTL}_{d,n}(u)$ if and only if one of the following cases holds:
\begin{enumerate}
\item [(i)] 
For some  $0 \leq m_1 \leq d-1$ we have $x_\ell=\exp(\ell m _1)$ ($0 \leq \ell \leq d-1$). In this case  the $x_i$'s are $d^{th}$ roots of unity and 
$z=-\frac{1}{u+1}$ or $z=-1$.
\item [(ii)] 
For some $0 \leq m_1,m_2 \leq d-1$ the $x_\ell$ are expressed as $x_\ell=\frac{1}{2}\left(\exp(\ell m _1) +\exp(\ell m _2) \right)$ ($0 \leq \ell \leq d-1$). In this case we have that  $z=-\frac{1}{2}$.
 \end{enumerate}
 
We note that in both cases the $x_i$'s are solutions of the ${\rm E}$--system, as required by \cite{JL2}, in order to proceed with defining link invariants. We do not take into consideration the case where: $z=-1$ (and the $x_i$'s are $d^{th}$ roots of unity) and the case where $z= - \frac{1}{2}$ (and $x_\ell=\frac{1}{2}\left(\exp(\ell m_1) + \exp(\ell m_2)  \right )
)$ since crucial braiding information is lost and therefore they are of no topological interest. Indeed, the trace ${\rm tr}$, for these two values of $z$ gives the same result for all even (resp. odd) powers of the $g_i$'s, as it becomes clear from the following formulas from \cite{JL2}, for $m \in \mathbb{Z}^{>0}$:
\begin{equation}\label{evenpower}
{\rm tr}(g_i^m) = \left( \frac{u^m -1}{u+1} \right )z + \left( \frac{u^m -1}{u+1} \right )E +1 \qquad \text{if } m \text{ is even}
\end{equation} 
and
\begin{equation}\label{oddpower}
{\rm tr}(g_i^m) = \left( \frac{u^m -1}{u+1} \right ) z + \left( \frac{u^m -1}{u+1} \right ) E - E \qquad \text{if } m \text{ is odd.}
\end{equation}
Notice that, substituting in Eq.~\ref{trg120} $z=-1$ implies $E=1$, while substituting $z=-\frac{1}{2}$ implies $E=\frac{1}{2}$.\\

The only remaining case of interest is case (i) where the $x_i$'s are roots of unity and $z= - \frac{1}{u+1}$. 
This implies that $E=1$ and $w=u$ in both Eqs.~\ref{gammainv} and \ref{deltainv}. We give the following definition:
\begin{definition}\label{invsdef} For $x_i$'s $d^{th}$ roots of unity ($x_i= x_1^i$, $ 1 \leq i \leq d-1$) and $z = -\frac{1}{u+1}$, we obtain from $\Gamma_{d,S}$ the following polynomial for $\alpha \in \cup_{\infty}\mathcal{F}_n$:
$$
\begin{array}{crcl}
 (i)&  \mathcal{V}_{d,S}(u)(\hat{\alpha})& =& \left (- \frac{1+ u}{\sqrt{u}}\right)^{n-1} \left(\sqrt{u}\right)^{\varepsilon(\alpha)} {\rm tr}(\gamma(\alpha)) = \Gamma_{d,S}(u,u). 
\end{array}
$$
Further, from $\Delta_{d,S}$, we obtain the following polynomial for $\alpha \in \cup_{\infty} B_n$:
$$
\begin{array}{crcl}
(ii) & V_{d,S}(u)(\hat{\alpha})& =& \left(- \frac{1+u}{\sqrt{u}} \right)^{n-1} \left(\sqrt{u}\right)^{\varepsilon(\alpha)} {\rm tr}(\delta(\alpha))= \Delta_{d,S}(u,u).
\end{array}
$$
Both polynomilas lie in $K(z,x_1,\ldots , x_{d-1})$.
\end{definition}
By Theorem~\ref{ytlthmgen} and the results of \cite{JL2} and \cite{JL4}, the polynomials $\mathcal{V}_{d,S}(u)$ and $V(u)$ are invariants of framed links and classical links respectively.\\

We know from \cite[Remark 5]{JL2} that the invariant $\Gamma_{d,S}(w,u)$ is not very interesting for framed links when the $x_i$'s are $d^{th}$ roots of unity because basic pairs of framed links are not distinguished. For classical links, as mentioned earlier, we know from \cite[Corollary 1]{ChLa} that the invariants $\Delta_{d,S}(w,u)$ coincide with the HOMFLYPT polynomial (case $E=1$). More precisely, for $E=1$ an algebra homomorphism can be defined, $h : {\rm Y}_{d,n}(u) \longrightarrow {\rm H}_n(u)$, and the composition $\tau \circ h$ is a Markov trace on ${\rm Y}_{d,n}(u)$ which takes the same values as the specialized trace ${\rm tr}$, whereby the $x_i$'s are specialized to $d^{th}$ roots of unity ($x_i= x_1^i$, $ 1 \leq m \leq d-1$). For details see \cite[\S3]{ChLa}.
 The above discussion leads to the following corollary:
 \begin{corollary}\label{invcor}
 The invariants $V_{d,S}(u)$ coincide with the Jones polynomial. The invariants $\mathcal{V}_{d,S}(u)$ are analogues of the Jones polynomial in the framed category.
 \end{corollary}

\end{document}